\documentclass{article}
\usepackage{amsmath, amsfonts, amssymb, amsthm}
\usepackage{graphicx, subcaption, cite}
\newtheorem{theorem}{Theorem}[section]
\newtheorem{lemma}[theorem]{Lemma}

\newtheorem{corollary}[theorem]{Corollary}
\newtheorem{exmp}[theorem]{Example}

\newcommand\bdf{\left\{\begin{array}{cccccc}}
\newcommand\edf{\end{array}\right.}
\newcommand\ba{\left(\begin{array}{cccccc}}

\newcommand\ea{\end{array}\right)}
\newcommand\bd{\left|\begin{array}{cccccc}}
\newcommand\ed{\end{array}\right|}
\newcommand\be{\begin{equation}}
\newcommand\ee{\end{equation}}
\newcommand\bes{\begin{equation}\begin{aligned}}
\newcommand\ees{\end{aligned}\end{equation}}
\newcommand\non{\nonumber}

\newcommand{\eq}[1]{Eq. \eqref{#1}}
\newcommand{\eqs}[2]{Eqs. (\ref{#1}) and (\ref{#2})}

\newcommand{\fig}[1]{Fig. \ref{#1}}
\newcommand{\tab}[1]{Table \ref{#1}}
\newcommand{\lem}[1]{Lemma \ref{#1}}
\newcommand{\thm}[1]{Theorem \ref{#1}}

\newcommand{\exm}[1]{Example \ref{#1}}

\newcommand{\half}{\frac{1}{2}}
\newcommand{\hf}[1]{\frac{#1}{2}}
\newcommand{\inv}[1]{\frac{1}{#1}}

\newcommand\diff[2]{\frac{\textrm{d}{#1}}{\textrm{d}{#2}}}
\newcommand\eps{\epsilon}
\newcommand\comment[1]{{}}

\title{Superconvergence of discontinuous Galerkin method for scalar and vector linear advection equations}
\author{Sirvan Rahmati, Tianshi Lu}

\begin{document}
\maketitle

\abstract{
In this paper, we use Fourier analysis to study the superconvergence of the semi-discrete discontinuous Galerkin method for scalar linear advection equations in one spatial dimension. The error bounds and asymptotic errors are derived for initial discretization by $L_2$ projection, Gauss-Radau projection, and other projections proposed by Cao et. al. \cite{Cao}. For pedagogical purpose, the errors are computed in two different ways. In the first approach, we compute the difference between the numerical solution and a special interpolation of the exact solution, and show that it consists of an asymptotic error of order $2k+1$ and a transient error of lower order. In the second approach, as in Ref. \cite{Cha}, we compute the error directly by decomposition into physical and nonphysical modes, and obtain agreement with the first approach. We then extend the analysis to vector conservation laws, solved using the Lax-Friedrichs flux. We prove that the superconvergence holds with the same order. The error bounds and asymptotic errors are demonstrated by various numerical experiments for scalar and vector advection equations.}

\section{Introduction}
Discontinuous Galerkin (DG) method is a class of finite element methods that uses discontinuous piecewise polynomials of order up to $k$ as test functions. The DG scheme is used widely for solving linear and nonlinear partial differential equations. Reed and Hill \cite{Reed} introduced the DG method for solving a steady-state linear hyperbolic equation in 1973. Cockburn et al. \cite{Cockburn1, Cockburn2, Cockburn3, Cockburn4} applied it to time-dependent nonlinear conservation laws.

In the past two decades, various superconvergence properties of DG methods have been studied, which provided a deeper understanding of DG solutions. According to the error in the DG method, the superconvergence divides into the following three categories. The first category is superconvergence of the error in cell average and at Radau points, measured in the discrete $L_2$ norm (see, e.g., \cite{Adjerid1,Adjerid2,Adjerid3}). The second category is the superconvergence of the DG solution towards a particular projection of the exact solution, called the supercloseness, typically measured in the standard $L_2$ norm. Some of the results are available in \cite{Cao1,Cheng,Yang,Baccouch,Meng}. The last category is the superconvergence of the post-processed solution. Negative-order norm estimates are a standard tool to derive superconvergent error estimates of the post-processed solution in the $L_2$ norm.  The choice of negative-order norms is to detect the oscillatory behavior of a function around zero. Post-processing aims to obtain a better approximation by convolving the numerical solution by a local averaging operator. For more information, we refer the reader to \cite{Bramble,Cockburn5,Mirzaee,Ji,Ji1}.
	
Now we shall review some of the superconvergence results and some of the relevant methods used in our work. Adjerid et al. \cite{Adjarid}  was first to show that the DG solution is superconvergent at Radau points for solving ordinary differential equations and steady-state hyperbolic problems. Later, Yang and Shu in \cite{Yang} studied superconvergence properties of the DG method for linear hyperbolic equations and proved that with suitable initial discretization, the error between the DG solution and the exact solution is $(k+2)$th order superconvergent at the downwind-biased Radau points. Around the same time,  Guo, Zhong, Qiu \cite{Guo} used the Fourier approach and decomposed the error. They symbolically computed eigenvalues and the corresponding eigenvectors of the DG method for low order approximations. Shortly after, Cao, Zhang, and Zou \cite{Cao1} showed that if the initial discretization is close enough to a particular reconstructed function, then the $2k+1$th (or $2k+1/2$th) superconvergence rate at the downwind points as well as the domain average, is achieved. Most studies have concentrated on the order of accuracy and neglected the role that the error coefficient plays in the estimates. Recently Frean and Ryan \cite{Frean} used a similar approach and showed that the semi-discrete error has dissipation errors of order $2k + 1$ and $2k + 2$ order of dispersion. They showed the critical role of the error constant in the dispersion and dissipation error for approximation polynomial degree $k$,  where $k = 0, 1, 2, 3$.
	
Cao showed the ($2k+1$) order of convergence in \cite{Cao1} for the first time, using the correction function technique. The idea of this technique is to construct a suitable correction function to correct the error between the exact solution and its Radau projection. Many papers have used this technique to show the semi-discrete DG method's superconvergence for one-dimensional problems. One of the recent papers done by Xu, Meng, Shu, Zhang \cite{Xu} uses a slightly modified correction function and the $L_2$-norm stability to establish the superconvergence property of the Runge-Kutta discontinuous Galerkin method for solving a linear constant-coefficient hyperbolic equation. They show that under a $r+1$ temporal and $2k+2$ spatial smoothness assumption, and by choosing a specific initialization, the cell average and the numerical flux are $\min(2k+1,r)$ superconvergent. They also prove a similar result for the post-processed solution, even if the initialization is the $L_2$ or Gauss-Radau projection.
	
Besides the Fourier analysis and the correction function technique, Pade approximation is another standard method for analyzing the DG method's superconvergence. Krivodonova and Qiu \cite{Kri} used Pade approximation to analyze the spectrum of the DG method when applied to the advection equation. They showed that for a uniform computational mesh of $ N $ elements, the eigenvalues could be classified into $ N $ physical modes and $ kN $ non-physical modes. They also show a $2k+2$ order of accuracy for the physical eigenvalue approximation. Later Chalmers and Krivodonova \cite{Cha} used Pade approximation to show the $(k+2)$th rate of superconvergence at the downwind points. They also proved that the $L_2$ projection of the numerical solution onto the $n$th Legendre polynomial is $2k+1-n$ accurate under certain initial conditions.

Most studies of the superconvergence for linear equations have been on single Fourier mode. The unphysical modes decay exponentially in time for each Fourier mode. However, the decay is exceedingly slow for high frequency modes. We analyze all modes as a whole, and prove the superconvergence in spite of the slow decay for high frequency modes. Fourier analysis and the correction function technique have been widely used to study the superconvergence of the numerical errors, yet the asymptotic error has not been given explicitly. The derivation of the asymptotic error can show the connection between Fourier analysis and the correction function technique, and clarify the effect of the initial discretization. The research on superconvergence has been focused on scalar equations. We study the superconvergence of the DG method applied to vector linear advection equations. To avoid diagonalization, the Lax-Friedrichs flux is used as the numerical flux on cell boundaries. The error analysis is adapted for the method since it is not equivalent to the upwinding flux for each scalar wave. We obtain similar results on the superconvergence.

The remainder of this paper is organized as follows. In Section 2, we study the superconvergence of the DG method for scalar linear advection equations. We use Fourier analysis to derive the error bounds. Then we compute the asymptotic errors for various initial discretization in two different approaches. In Section 3, we compute the numerical error for the DG method with Lax-Friedrichs flux applied to vector advection equations. We show the superconvergence similar to that for scalar equations. The results are illustrated numerically in Section 4. Finally, conclusions and thoughts on future work are discussed in Section 5.

\section{Scalar Equations}
In this section we investigate the superconvergence of the Discontinuous Galerkin (DG) method for the one-dimensional scalar linear advection equation with periodic boundary condition.

\subsection{Preliminary}
The linear advection equation is
\bes\label{master}
&u_t+au_x=0,\ \ \ (x,t)\in[0,2\pi]\times[0,T].\\
&u(x,0)=g(x),\ \ \ u(0,t)=u(2\pi,t).
\ees
Without loss of generality, the sound speed $a$ is set to 1. We solve the equation using $k$'th order DG method on a uniform grid of $N$ cells with cell size $h$, and denote the numerical solution by $u_h$, The solution $u_h\in V_h$, where
\be
V_h=\{v: v|_{\tau_j}\in P_k(\tau_j), 1\le j\le N\},
\ee
and satisfies
\be\label{numeq}
((u_h)_t,v)=(u_h,v_x)+\sum_{j=1}^N([v]u_h^-)|_{j+\half},\ \ \ \forall v\in V_h,
\ee
where $[v]=v^+-v^-$, and $v^+_{N+\half}\equiv v^+_{\half}$. The initial discretization $u_h(x,0)$ is a projection of $g(x)$ onto $V_h$. To analyze the error $e=u-u_h$, we follow Cao \cite{Cao} to define the special interpolations of the exact solution,
\be\label{def-uI}
u_I=P^-_hu-\sum_{i=1}^kw_i,
\ee
where $P^-_h$ is the Gauss-Radau projection, and in each cell $\tau_j$, for $1\le i\le k$,
\bes\label{def-w}
(w_i,v_x)&=(\partial_tw_{i-1},v),\ \ \ \forall v\in P_k(\tau_j),\\
w_i^-&=0,
\ees
with $w_0=u-P^-_hu$. Let $\eta=u-u_I$, and $e=\xi+\eta$. The equation for $\xi$ is
\be\label{eq-xi}
((\xi+w_k)_t,v)=(\xi,v_x)+\sum_{j=1}^N([v]\xi^-)|_{j+\half},\ \ \ \forall v\in V_h,
\ee
We can express $e$, $\eta$, and $\xi$ in terms of orthogonal basis $\{\phi_{j,n}\}$, Legendre polynomials \cite{Abr} scaled to $\tau_j$, as
\be
e|_{\tau_j}(t)=\sum_{n=0}^\infty e_{j,n}(t)\phi_{j,n},\ \ \
\eta|_{\tau_j}=\sum_{n=0}^\infty\eta_{j,n}\phi_{j,n},\ \ \
\xi|_{\tau_j}(t)=\sum_{n=0}^k \xi_{j,n}(t)\phi_{j,n}.
\ee
It has been proven by Cao \cite{Cao} that for $t\in[0,T]$ and $0\le n\le k$,
\be\label{eta-bd}
\|\eta_{\cdot,n}(t)\|_\infty\lesssim h^{2k+1-n}\|g^{(2k+1-n)}\|_\infty,
\ee
provided that $g^{(2k+1)}(x)$ is bounded. Cao \cite{Cao} also showed that if the initialization is chosen such that
\be
\|u_h(\cdot,0)-u_I(\cdot,0)\|_2\lesssim h^{2k+1}\|g^{(2k+2)}\|_\infty,
\ee
then
\be
\|\xi(\cdot,t)\|_2\lesssim h^{2k+1}\|g^{(2k+2)}\|_\infty,
\ee
for $t\in[0,T]$, provided that $g^{(2k+2)}(x)$ is bounded. Consequently, for $t\in[0,T]$ and $0\le n\le k$,
\be
\|e_{\cdot,n}(t)\|_2\lesssim h^{2k+1-n}\|g^{(2k+2)}\|_\infty.
\ee
In addition, the downwind error,
\be
e^-_j(t)=u(x^-_{j+\half},t)-u_h(x^-_{j+\half},t),
\ee
is also of order $2k+1$,
\be
\|e^-_\cdot(t)\|_2\lesssim h^{2k+1}\|g^{(2k+2)}\|_\infty.
\ee
We will apply Fourier analysis to \eq{master}, and prove the we can obtain optimal superconvergence asymptotically by initializing $u_h$ as the $L_2$ projection of $g$ onto $P_k(\tau_j)$ on each cell.

\subsection{Error estimation for initialization by $L_2$ projection}

For a discrete function $f_j$, we denote the Fourier coefficients by $\hat{f}_m$,
\be
\hat{f}_m=\inv{N}\sum_{j=1}^Nf_je^{-imx_j}.
\ee
and define the norm,
\be
\|f\|_s=\sum_{m=0}^{N-1}|\hat{f}_m|m^s.
\ee
The next theorem bounds the error in $u_h(x,t)$ for $u_h(x,0)$ close to the $L_2$ projection of $g(x)$.

\begin{theorem}\label{Thm-ebd}
Suppose $g\in H_{2k+3}$, so that $\|g\|_{2k+2}<\infty$. Let $u_h$ be the solution to \eq{numeq} with initial error
\be\label{einit}
\hat{e}_{m,n}(0)=O((mh)^{2k+1-n}\hat{g}_m),\ \ \ 0\le n\le k.
\ee
There exists an $\alpha>0$, such that for any $t\in[0,T]$,
\be\label{xi-order}
\|\xi(\cdot,t)\|_\infty\lesssim h^{2k+1}(\|g\|_{2k+1}+\|g\|_{2k+2}t)+h^{k+1}\|g\|_{k+1}
e^{-\frac{\alpha t}{2h}}.
\ee
\be\label{xi0-order}
\|\xi_{\cdot,0}(t)\|_\infty\lesssim h^{2k+1}(\|g\|_{2k+1}+\|g\|_{2k+2}t)+h^{k+2}\|g\|_{k+2}
e^{-\frac{\alpha t}{2h}}.
\ee
Consequently,
\be\label{e0-order}
\|e_{\cdot,0}(t)\|_\infty\lesssim h^{2k+1}(\|g\|_{2k+1}+\|g\|_{2k+2}t)+h^{k+2}\|g\|_{k+2}
e^{-\frac{\alpha t}{2h}},
\ee
\be\label{er-order}
\|e^-_\cdot(t)\|_\infty\lesssim h^{2k+1}(\|g\|_{2k+1}+\|g\|_{2k+2}t)+h^{k+1}\|g\|_{k+1}
e^{-\frac{\alpha t}{2h}},
\ee
and for $1\le n\le k$,
\be\label{en-order}
\|e_{\cdot,n}(t)\|_\infty\lesssim h^{2k+1-n}\|g\|_{2k+1-n}+h^{k+1}\|g\|_{k+1}
e^{-\frac{\alpha t}{2h}}.
\ee
\end{theorem}

\begin{proof}
For each Fourier mode, \eq{eq-xi} can be written as a matrix equation,
\be
\diff{}{t}\hat{\xi}_m=\frac{A_m}{h}\hat{\xi}_m-b_m(t),
\ee
where $\hat{\xi}_m=(\hat{\xi}_{m,0},\hat{\xi}_{m,1},\ldots,\hat{\xi}_{m,k})^T$,
\be
b_m(t)=\diff{}{t}((\hat{w}_k)_{m,0},(\hat{w}_k)_{m,1},\ldots,(\hat{w}_k)_{m,k})^T,
\ee
and $A_m$ is a $(k+1)\times(k+1)$ matrix with indices from 0 to $k$,
\bes
(A_m)_{st}&=-(2s+1)(1-(-1)^s e^{-imh}),\ \ \ 0\le s\le t\le k,\\
(A_m)_{st}&=-(2s+1)(-1)^{s+t}(1-(-1)^t e^{-imh}),\ \ \ 0\le t<s\le k.
\ees
The solution is
\be\label{xi-sol}
\hat{\xi}_m(t)=e^{A_mt}\hat{\xi}_m(0)-\int_0^te^{A_m(t-\tau)}b_m(\tau)d\tau.
\ee
Diagonalize $A_m$ as $A_m=\sum_{n=0}^k\lambda_n^{(m)}r_n^{(m)}l_n^{(m)}$, where $\lambda_n^{(m)}$, $r_n^{(m)}$, and $l_n^{(m)}$ are the eigenvalues and the associated right and left eigenvectors that satisfy $l_n^{(m)}r_n^{(m)}=1$. Because
\be
\frac{(-1)^s(A_m)_{st}}{2s+1}=\frac{(-1)^t(A_m)_{ts}}{2t+1},
\ee
we can set
\be\label{lr}
(l_n^{(m)})_s=\frac{(-1)^s}{2s+1}(r_n^{(m)})_s,\ \ \ 0\le s\le k.
\ee
It have been shown in Ref. \cite{Cha} that $\lambda_0$ represents the physical mode,
\be\label{l0}
\lambda^{(m)}_0=-imh+O((mh)^{2k+2}),
\ee
while other eigenvalues ($1\le n\le k$) represent nonphysical modes,
\be\label{nonphy}
\lambda^{(m)}_n=-\alpha_n+O(mh),\ \ \ \Re\alpha_n>0.
\ee
In fact, $\Re\lambda^{(m)}_n<0$ for $1\le n\le k$. To prove that, notice that \eq{numeq} can be written as
\be\label{intpart}
((u_h)_t+(u_h)_x,v)=-\sum_{j=1}^N([u_h]v^+)|_{j+\half},\ \ \ \forall v\in V_h.
\ee
It implies the energy estimate,
\be\label{energy}
\diff{}{t}(u_h,\bar{u_h})=-\sum_{j=1}^N|[u_h]_{j+\half}|^2\le0.
\ee
Consequently, $\Re\lambda^{(m)}_n\le0$. If $\Re\lambda^{(m)}_n=0$, $[u_h]=0$ by \eq{energy}. Then we have $(u_h)_t+(u_h)_x=0$ by \eq{intpart}. As an eigenfunction in $P_k$, $u_h$ has to be a constant, and the associated eigenvalue is 0. But that is the physical eigenfunction for $m=0$, so nonphysical eigenvalues have negative real parts. We can write
\be\label{xi1}
e^{A_mt}\hat{\xi}_m(0)=e^{\lambda_0^{(m)}t/h}r_0^{(m)}l_0^{(m)}\hat{\xi}_m(0)+\sum_{n=1}^k e^{\lambda_n^{(m)}t/h}r_n^{(m)}l_n^{(m)}\hat{\xi}_m(0).
\ee
Express $e^{im(x-x_j)}$ in terms of orthogonal basis $\{\phi_{j,n}\}$,
\be
e^{im(x-x_j)}=\sum_{n=0}^\infty p_n^{(m)}\phi_{j,n}(x).
\ee
For the physical mode, it has been shown in Ref. \cite{Cha} that
\be\label{r0}
(r_0^{(m)})_n=p_n^{(m)}+O((mh)^{2k+1-n})=O((mh)^n),\ \ \ 0\le n\le k,
\ee
where we used
\be\label{pn}
p_n^{(m)}=\frac{n!}{(2n)!}(imh)^n+O((mh)^{n+1}).
\ee
Then $(l_0^{(m)})_n=O((mh)^n)$, and
\be
l_0^{(m)}\hat{\xi}_m(0)=\sum_{n=0}^k(\hat{e}_{m,n}(0)-\hat{\eta}_{m,n}(0))O((mh)^n).
\ee
By \eq{eta-bd}, $\hat{\eta}_{m,n}(0)=O((mh)^{2k+1-n}\hat{g}_m)$. Combined with \eq{einit}, we have
\be\label{xi-phy}
e^{\lambda_0^{(m)}t/h}r_0^{(m)}l_0^{(m)}\hat{\xi}_m(0)=O((mh)^{2k+1}\hat{g}_m).
\ee
For nonphysical modes, by \eq{nonphy} there exists $\eps>0$ such that for $mh\le\eps$,
\be
\Re\lambda_n\le-\hf{\alpha},\ \ \ 1\le n\le k,
\ee
where
\be
\alpha=\min_{1\le n\le k}\Re\alpha_n.
\ee
For $mh\le\eps$,
\be\label{xi-nonphy}
e^{\lambda_n^{(m)}t/h}r_n^{(m)}l_n^{(m)}\hat{\xi}_m(0)
=r_n^{(m)}O(e^{-\frac{\alpha t}{2h}}(mh)^{k+1}\hat{g}_m).
\ee
whose first entry is
\be\label{xi10}
\left[e^{\lambda_n^{(m)}t/h}r_n^{(m)}l_n^{(m)}\hat{\xi}_m(0)\right]_0
=O(e^{-\frac{\alpha t}{2h}}(mh)^{k+2}\hat{g}_m).
\ee
because $(r_n^{(m)})_0=O(mh)$. For $mh>\eps$,
\be
e^{\lambda_n^{(m)}t/h}r_n^{(m)}l_n^{(m)}\hat{\xi}_m(0)=O((mh)^{k+1}\hat{g}_m).
\ee
For the second term in \eq{xi-sol}, by the definition of $w_k$ in \eq{def-w}, we have
\be\label{bm}
b_m=O(m^{2k+2}h^{2k+1}\hat{g}_m),
\ee
and so
\be\label{xi2}
\int_0^te^{A_m(t-\tau)}b_m(\tau)d\tau=O(m^{2k+2}h^{2k+1}\hat{g}_mt),
\ee
Summing over all Fourier modes, we obtain
\be
\|\xi(\cdot,t)\|_\infty \le\sum_{n=0}^k\|\xi_{\cdot,n}(t)\|_\infty
\le\sum_{n=0}^k\sum_{m=0}^{N-1}|\hat{\xi}_{m,n}|
\le\sum_{n=0}^k\sum_{m=0}^{N-1}
|[e^{A_mt}\hat{\xi}_m(0)]_n|+|[\int_0^te^{A_m(t-\tau)}b_m(\tau)d\tau]_n|.
\ee
Since for any $0\le s\le k$,
\bes\label{xi1-bd}
\left|\sum_{m=0}^{N-1}[e^{A_mt}\hat{\xi}_m(0)]_s\right|&\lesssim \sum_{m=0}^{N-1}(mh)^{2k+1}|\hat{g}_m|+
\sum_{m=0}^{\eps/h}|(r_n)_s|(mh)^{k+1}|\hat{g}_m|e^{-\frac{\alpha t}{2h}}+
\sum_{m=\eps/h}^{N-1}(mh)^{k+1}|\hat{g}_m|\\
&\le \sum_{m=0}^{N-1}(mh)^{2k+1}|\hat{g}_m|+
\sum_{m=0}^{\eps/h}|(r_n)_s|(mh)^{k+1}|\hat{g}_m|e^{-\frac{\alpha t}{2h}}+
\eps^{-k}\sum_{m=\eps/h}^{N-1}(mh)^{2k+1}|\hat{g}_m|\\
&\lesssim h^{2k+1}\|g\|_{2k+1}+h^{k+1}\|g\|_{k+1}e^{-\frac{\alpha t}{2h}}.
\ees
and
\be\label{xi2-bd}
\sum_{m=0}^{N-1}|[\int_0^te^{A_m(t-\tau)}b_m(\tau)d\tau]_n|=O(h^{2k+1}\|g\|_{2k+2}t).
\ee
we obtain \eq{xi-order}. For $s=0$, by \eq{xi10},
\be\label{xi10-bd}
\sum_{m=0}^{N-1}[e^{A_mt}\hat{\xi}_m(0)]_0
=O(h^{2k+1}\|g\|_{2k+1}+h^{k+2}\|g\|_{k+2}e^{-\frac{\alpha t}{2h}}),
\ee
hence \eq{xi0-order}. \eq{e0-order} is a consequence of \eqs{eta-bd}{xi0-order}. Since $\eta_j^-=0$,
\be
e_j^-(t)=\xi_j^-(t)=\sum_{n=0}^k\xi_{j,n}(t),
\ee
hence \eq{er-order}. \eq{en-order} is a consequence of \eqs{eta-bd}{xi-order}.
\end{proof}

If we replace the norm $\|f\|_s$ by $\|f^{(s)}\|_{L_2}$, or equivalently,
\be
\|f\|_{s,2}\equiv\sqrt{\sum_{m=0}^{N-1}|\hat{f}_mm^s|^2},
\ee
the proof of \thm{Thm-ebd} can be carried out similarly, which gives the following error bound in $L_2$ norm.

\begin{theorem}\label{Thm-ebd2}
Suppose $g\in H_{2k+2}$, so that $\|g^{(2k+2)}\|_{L_2}<\infty$. Let $u_h$ be the solution to \eq{numeq} with initial error
\be
\hat{e}_{m,n}(0)=O((mh)^{2k+1-n}\hat{g}_m),\ \ \ 0\le n\le k.
\ee
There exists an $\alpha>0$, such that for any $t\in[0,T]$,
\be
\|\xi(\cdot,t)\|_{L_2}\lesssim h^{2k+1}
(\|g^{(2k+1)}\|_{L_2}+\|g^{(2k+2)}\|_{L_2}t)+h^{k+1}\|g^{(k+1)}\|_{L_2}e^{-\frac{\alpha t}{2h}}.
\ee
\be
\|\xi_{\cdot,0}(t)\|_{L_2}\lesssim h^{2k+1}
(\|g^{(2k+1)}\|_{L_2}+\|g^{(2k+2)}\|_{L_2}t)+h^{k+2}\|g^{(k+2)}\|_{L_2}e^{-\frac{\alpha t}{2h}}.
\ee
Consequently,
\be
\|e_{\cdot,0}(t)\|_{L_2}\lesssim h^{2k+1}
(\|g^{(2k+1)}\|_{L_2}+\|g^{(2k+2)}\|_{L_2}t)+h^{k+2}\|g^{(k+2)}\|_{L_2}e^{-\frac{\alpha t}{2h}},
\ee
\be
\|e^-_\cdot(t)\|_{L_2}\lesssim h^{2k+1}
(\|g^{(2k+1)}\|_{L_2}+\|g^{(2k+2)}\|_{L_2}t)+h^{k+1}\|g^{(k+1)}\|_{L_2}e^{-\frac{\alpha t}{2h}},
\ee
and for $1\le n\le k$,
\be
\|e_{\cdot,n}(t)\|_{L_2}\lesssim h^{2k+1-n}
\|g^{(2k+1-n)}\|_{L_2}+h^{k+1}\|g^{(k+1)}\|_{L_2}e^{-\frac{\alpha t}{2h}}.
\ee
\end{theorem}

\subsection{Asymptotic error}

We derive the asymptotic error as $h\to0$ for sufficiently small initial error.

\begin{theorem}\label{Thm-elim}
Let $u_h$ be the solution to \eq{numeq} with initial error
\be\label{einit2}
\hat{e}_{m,n}(0)=O((mh)^{2k+2-n}\hat{g}_m),\ \ \ 0\le n\le k.
\ee
For any $t>0$,
\be\label{e0-lim}
\lim_{h\to0}\frac{e_{j,0}(t)}{h^{2k+1}}=\frac{(-1)^k(k+1)!k!}{(2k+2)!(2k+1)!}
[kg^{(2k+1)}(x_j-t)-tg^{(2k+2)}(x_j-t)].
\ee
\be\label{er-lim}
\lim_{h\to0}\frac{e^-_j(t)}{h^{2k+1}}=\frac{(-1)^k(k+1)!k!}{(2k+2)!(2k+1)!}
[kg^{(2k+1)}(x_j-t)-tg^{(2k+2)}(x_j-t)].
\ee
The convergence is uniform if $g\in H_{2k+3}$, or in $L_2$ norm if $g\in H_{2k+2}$. For $1\le n\le k$,
\be\label{en-lim}
\lim_{h\to0}\frac{e_{j,n}(t)}{h^{2k+1-n}}=
\frac{(-1)^{k+1-n}(k+1)!k!(2n+1)!}{(2k+2)!(2k+1)!n!}g^{(2k+1-n)}(x_j-t).
\ee
The convergence is uniform if $g\in H_{2k+2}$.
\end{theorem}

\begin{proof}
First we assume $g\in H_{2k+3}$, so that $\|g\|_{2k+2}<\infty$. By \eq{xi-sol}, the error vector in each cell is
\be\label{e-sol}
e_j(t)=\sum_{m=0}^{N-1}e^{imx_j}\left[e^{A_mt}\hat{\xi}_m(0)
-\int_0^te^{A_m(t-\tau)}b_m(\tau)d\tau+\hat{\eta}_m(t)\right].
\ee
By the construction of $u_I$ in \eq{def-uI}, for smooth $g$,
\be
\eta_{j,k}=(w_0)_{j,k}+O(h^{k+2}\|g^{(k+2)}\|_\infty)
=-\frac{(k+1)!}{(2k+2)!}h^{k+1}\partial_x^{k+1}u+O(h^{k+2}\|g^{(k+2)}\|_\infty).
\ee
Similarly, for $0\le n\le k-1$,
\bes
\eta_{j,n}&=(w_{k-n})_{j,n}+O(h^{2k+2-n}\|g^{(2k+2-n)}\|_\infty)\\
&=-\frac{(k+1)!k!(2n+1)!}{(2k+2)!(2k+1)!n!}
h^{2k+1-n}\partial_t^{k-n}\partial_x^{k+1}u+O(h^{2k+2-n}\|g^{(2k+2-n)}\|_\infty).
\ees
Therefore, for each Fourier mode and $0\le n\le k$,
\be\label{eta-t}
\hat{\eta}_{m,n}(t)=-\frac{(k+1)!k!(2n+1)!}{(2k+2)!(2k+1)!n!}
h^{2k+1-n}(-im)^{k-n}(im)^{k+1}\hat{g}_me^{-imt}+O((mh)^{2k+2-n}\hat{g}_m).
\ee
By \eq{einit2},
\be\label{xi0}
\hat{\xi}_{m,n}(0)=\frac{(k+1)!k!(2n+1)!}{(2k+2)!(2k+1)!n!}
h^{2k+1-n}(-im)^{k-n}(im)^{k+1}\hat{g}_m+O((mh)^{2k+2-n}\hat{g}_m).
\ee
By \eq{xi1},
\be
e^{A_mt}\hat{\xi}_m(0)=e^{\lambda_0^{(m)}t/h}r_0^{(m)}l_0^{(m)}\hat{\xi}_m(0)+\sum_{n=1}^k e^{\lambda_n^{(m)}t/h}r_n^{(m)}l_n^{(m)}\hat{\xi}_m(0).
\ee
By \eqs{lr}{pn},
\be
(l_0^{(m)})_n=\frac{(-1)^n}{2n+1}p_n^{(m)}+O((mh)^{2k+1-n})=\frac{(-1)^nn!}{(2n+1)!}(imh)^n+O((mh)^{n+1}).
\ee
So
\be\label{l0xi}
l_0^{(m)}\hat{\xi}_m(0)=(k+1)\frac{(k+1)!k!}{(2k+2)!(2k+1)!}
(mh)^{2k+1}i\hat{g}_m+O((mh)^{2k+2}\hat{g}_m).
\ee
Since $(r_0^{(m)})_0=1+O(mh)$, by \eq{l0},
\bes
&e^{\lambda_0^{(m)}t/h}(r_0^{(m)})_0l_0^{(m)}\hat{\xi}_m(0)=e^{-imt+O(mh)mt}
\left(\frac{(k+1)(k+1)!k!}{(2k+2)!(2k+1)!}(mh)^{2k+1}i\hat{g}_m+O((mh)^{2k+2}\hat{g}_m)\right)\\
&=e^{-imt}\frac{(k+1)(k+1)!k!}{(2k+2)!(2k+1)!}(mh)^{2k+1}i\hat{g}_m+O((mh)^{2k+2}\hat{g}_m)
+O((mh)^{2k+1}\hat{g}_m)(e^{O(mh)mt}-1).
\ees
The last two terms are $o(h^{2k+1})$ when summed over all Fourier modes, because
\be
\frac{\sum_{m=0}^{N-1}|(mh)^{2k+2}\hat{g}_m|}{h^{2k+1}}=\sum_{m=0}^{1/\sqrt{h}}m^{2k+1}|\hat{g}_m|mh
+\sum_{m=1/\sqrt{h}}^{N-1}m^{2k+1}|\hat{g}_m|mh\le\|g\|_{2k+1}\sqrt{h}+\sum_{m=1/\sqrt{h}}^\infty m^{2k+1}|\hat{g}_m|,
\ee
which converges to 0 as $h\to0$ since $\|g\|_{2k+1}<\infty$, and as $\Re\lambda_0\le0$,
\bes
\frac{\sum_{m=0}^{N-1}(mh)^{2k+1}|\hat{g}_m(e^{O(mh)mt}-1)|}{h^{2k+1}}
&\lesssim\sum_{m=0}^{1/\sqrt[3]{h}}m^{2k+1}|\hat{g}_m|m^2ht
+\sum_{m=1/\sqrt[3]{h}}^{N-1}m^{2k+1}|\hat{g}_m|\\
&\le\|g\|_{2k+1}\sqrt[3]{h}t+\sum_{m=1/\sqrt[3]{h}}^\infty m^{2k+1}|\hat{g}_m|,
\ees
which also converges to 0 as $h\to0$. For $1\le s\le k$, since $(r_0^{(m)})_s=O((mh)^s)$,
\be
\sum_{m=0}^{N-1}e^{\lambda_0^{(m)}t/h}(r_0^{(m)})_sl_0^{(m)}\hat{\xi}_m(0)=o(h^{2k+1}).
\ee
Similar to \eq{xi1-bd}, for $1\le n\le k$ and $0\le s\le k$, we have
\bes
\inv{h^{2k+1}}\left|\sum_{m=0}^{N-1}e^{\lambda_n^{(m)}t/h}(r_n^{(m)})_sl_n^{(m)}\hat{\xi}_m(0)\right|
&\lesssim \sum_{m=0}^{\eps/h}m^{k+1}h^{-k}|\hat{g}_m|e^{-\frac{\alpha t}{2h}}+
\sum_{m=\eps/h}^{N-1}m^{k+1}h^{-k}|\hat{g}_m|\\
&\le h^{-k}\|g\|_{k+1}e^{-\frac{\alpha t}{2h}}+
\eps^{-k}\sum_{m=\eps/h}^{N-1}m^{2k+1}|\hat{g}_m|,
\ees
which converges to 0 as $h\to0$. Therefore
\be\label{ei}
\sum_{m=0}^{N-1}e^{imx_j}[e^{A_mt}\hat{\xi}_m(0)]_0=
\sum_{m=0}^{N-1}e^{im(x_j-t)}\frac{(k+1)(k+1)!k!}{(2k+2)!(2k+1)!}(mh)^{2k+1}i\hat{g}_m+o(h^{2k+1}),
\ee
and $\sum_{m=0}^{N-1}e^{imx_j}[e^{A_mt}\hat{\xi}_m(0)]_s=o(h^{2k+1})$ for $1\le s\le k$. For the integral in \eq{e-sol}, we have \eq{bm}. In particular, since for smooth $g$,
\be
\diff{}{t}(w_k)_{j,k}=
-\frac{(k+1)!k!}{(2k+2)!(2k+1)!}h^{2k+1}\partial_t^{k+1}\partial_x^{k+1}u+O(h^{2k+2}\|g^{(2k+3)}\|_\infty).
\ee
we have
\be
(b_m)_0(\tau)=-\frac{(k+1)!k!}{(2k+2)!(2k+1)!}
h^{2k+1}(-im)^{k+1}(im)^{k+1}\hat{g}_me^{-im\tau}+O(m^{2k+3}h^{2k+2}\hat{g}_m).
\ee
Then
\bes
&e^{\lambda_0^{(m)}(t-\tau)/h}(r_0^{(m)})_0l_0^{(m)}b_m(\tau)\\
&=e^{-imt+O(mh)m(t-\tau)}
\left(-\frac{(k+1)!k!}{(2k+2)!(2k+1)!}(mh)^{2k+1}m\hat{g}_m+O((mh)^{2k+2}m\hat{g}_m)\right)\\
&=e^{-imt}\frac{-(k+1)!k!(mh)^{2k+1}}{(2k+2)!(2k+1)!}m\hat{g}_m+O((mh)^{2k+2}m\hat{g}_m)
+O((mh)^{2k+1}m\hat{g}_m)(e^{O(mh)mt}-1).
\ees
Similarly, the last two terms are $o(h^{2k+1})$ when summed over all Fourier modes because $\|g\|_{2k+2}<\infty$, and for $1\le s\le k$,
\be
\sum_{m=0}^{N-1}e^{\lambda_0^{(m)}(t-\tau)/h}(r_0^{(m)})_sl_0^{(m)}b_m(\tau)=o(h^{2k+1}).
\ee
For $1\le n\le k$,
\be
\inv{h^{2k+1}}\left|\sum_{m=0}^{N-1}e^{\lambda_n^{(m)}(t-\tau)/h}(r_n^{(m)})_0l_n^{(m)}b_m(\tau)\right|
\lesssim \sum_{m=0}^{N-1}m^{2k+3}h|\hat{g}_m|,
\ee
which converges to 0 at $h\to0$. Therefore
\be\label{es}
-\sum_{m=0}^{N-1}e^{imx_j}\int_0^t[e^{A_m(t-\tau)}b_m(\tau)]_0d\tau=
\sum_{m=0}^{N-1}e^{im(x_j-t)}\frac{(k+1)!k!}{(2k+2)!(2k+1)!}(mh)^{2k+1}mt\hat{g}_m+o(h^{2k+1}).
\ee
For $1\le s\le k$,
\bes
&\left|\sum_{m=0}^{N-1}e^{imx_j}\int_0^t[e^{A_m(t-\tau)}b_m(\tau)]_sd\tau\right|\\
&=\left|\sum_{n=1}^k\sum_{m=0}^{N-1}e^{imx_j}\int_0^t
e^{\lambda_n^{(m)}(t-\tau)/h}(r_n^{(m)})_sl_n^{(m)}b_m(\tau)d\tau+o(h^{2k+1})\right|\\
&\lesssim\sum_{n=1}^k\left[\sum_{m=0}^{\eps/h}\int_0^t
e^{-\frac{\alpha(t-\tau)}{2h}}m^{2k+2}h^{2k+1}|\hat{g}_m|d\tau+
\sum_{m=\eps/h}^{N-1}\int_0^tm^{2k+2}h^{2k+1}|\hat{g}_m|d\tau\right]+o(h^{2k+1})\\
&\le\sum_{n=1}^kh^{2k+1}\left[\frac{2h}{\alpha}\sum_{m=0}^{N-1}m^{2k+2}|\hat{g}_m|+
t\sum_{m=\eps/h}^{N-1}m^{2k+2}|\hat{g}_m|\right]+o(h^{2k+1})\\
&=o(h^{2k+1}).
\ees
Lastly, by \eq{eta-t},
\be\label{e3}
\hat{\eta}_{m,0}(t)=-\frac{(k+1)!k!}{(2k+2)!(2k+1)!}
(mh)^{2k+1}i\hat{g}_me^{-imt}+O((mh)^{2k+2}\hat{g}_m).
\ee
whose last term is $o(h^{2k+1})$ when summed over all Fourier modes as $\|g\|_{2k+1}<\infty$. Substituting \eq{ei}, \eqref{es} and \eqref{e3} into \eq{e-sol}, we obtain
\bes
\lim_{h\to0}\frac{e_{j,0}(t)}{h^{2k+1}}&=\lim_{N\to\infty}\sum_{m=0}^{N-1}\frac{(k+1)!k!}{(2k+2)!(2k+1)!}
[km^{2k+1}i+tm^{2k+2}]\hat{g}_me^{im(x_j-t)}\\
&=\frac{(-1)^k(k+1)!k!}{(2k+2)!(2k+1)!}[kg^{(2k+1)}(x_j-t)-tg^{(2k+2)}(x_j-t)].
\ees
For $g\in H_{2k+3}$, the Fourier series of $g^{(2k+2)}$ converges absolutely and uniformly. For $g\in H_{2k+2}$, the Fourier series of $g^{(2k+2)}$ converges in $L_2$ norm, and the proof is still valid so long as we replace all $\|g\|_s$ by $\|g\|_{s,2}$. To prove \eq{er-lim}, we notice that $\eta_j^-=0$, so
\be
e_j^-(t)=\xi_j^-(t)=\sum_{n=0}^k\xi_{j,n}(t)=\xi_{j,0}(t)+o(h^{2k+1}).
\ee
For $g\in H_{2k+2}$, we also have
\be
\|g\|_{2k+1}\lesssim\|g\|_{2k+2,2}<\infty.
\ee
Substituting \eq{bm} into \eq{e-sol} we get
\be
e_j(t)=O(h^{2k+1}\|g\|_{2k+1})+o(h^{2k}\|g\|_{2k+1})+\sum_{m=0}^{N-1}e^{imx_j}\hat{\eta}_m(t).
\ee
Substituting in \eq{eta-t} we obtain for $1\le n\le k$,
\bes
\lim_{h\to0}\frac{e_{j,n}(t)}{h^{2k+1-n}}&=\lim_{N\to\infty}\sum_{m=0}^{N-1}
\frac{(k+1)!k!(2n+1)!}{(2k+2)!(2k+1)!n!}
i^{n-1}m^{2k+1-n}\hat{g}_me^{im(x_j-t)}\\
&=\frac{(-1)^{k+1-n}(k+1)!k!(2n+1)!}{(2k+2)!(2k+1)!n!}g^{(2k+1-n)}(x_j-t),
\ees
with uniform convergence.
\end{proof}

\subsection{Direction computation of error}

We can provide an alternative proof of \thm{Thm-elim} by direct computation of the error. As shown in \eqs{l0}{r0}, the physical mode of $A_m$ is super close to the exact solution. Ref. \cite{Kri} proved that the $R_{k,k+1}(\lambda_0^{(m)})=\exp(imh)$, where $R_{k,k+1}(z)$ is the $[k/k+1]$ Pad\'{e} approximation of $\exp(-z)$. We will prove a lemma on the closeness of the physical mode to the exact solution following Ref. \cite{Cha}. Denote by $\tilde{r}_0^{(m)}$ the eigenvector parallel to $r_0^{(m)}$ but normalized at the downwind point,
\be\label{r0dw}
\sum_{n=0}^k(\tilde{r}_0^{(m)})_n=e^{\hf{imh}}.
\ee

\begin{lemma}\label{Lem-eig}
The eigenvalue of the physical mode is superclose to $-imh$. More precisely,
\be\label{lambda0}
\lambda_0^{(m)}=-imh-\frac{(k+1)!k!}{(2k+2)!(2k+1)!}(mh)^{2k+2}+O((mh)^{2k+3}).
\ee
The eigenvector normalized as in \eq{r0dw} is superclose to the projection of $e^{im(x-x_j)}$ onto $P_k(\tau_j)$. More precisely,
\be\label{delta}
\delta_n^{(m)}\equiv p_n^{(m)}-(\tilde{r}_0^{(m)})_n
=(-1)^{k+1-n}\frac{(k+1)!k!(2n+1)!}{(2k+2)!(2k+1)!n!}(imh)^{2k+1-n}+O((mh)^{2k+2-n})
\ee
for $0\le n\le k$.
\end{lemma}

\begin{proof}
Following Ref. \cite{Cha}, for each Fourier mode, we can write the eigenvalue problem for \eq{numeq}, scaled from $x\in\tau_j$ to $y=2(x-x_j)/h\in[-1,1]$, as
\be\label{eig-eq}
\lambda u+2u_y=(-1)^{k+1}[[u]](R_{k+1}^-)'(y),\ \ \ u\in P_k,\ u(1)=e^{\hf{imh}}.
\ee
where $[[u]]=u(-1)-u(1)e^{-imh}$, and $R^-_{k+1}(y)=\phi_{k+1}(y)-\phi_k(y)$ is the right Rado polynomial \cite{Abr} of degree $k+1$. The solution $u$ associated with the physical eigenvalue $\lambda_0^{(m)}$ is $\tilde{r}_0^{(m)}$. Since $u\in P_k$,
\be
u(y)=[[u]]\hf{(-1)^k}\sum_{l=1}^{k+1}\frac{R^{-,(l)}_{k+1}(y)}{(-\hf{\lambda})^l}.
\ee
Substituting in $y=1$ and using the formula $R^{-,(k+1)}_{k+1}=\phi^{(k+1)}_{k+1}=(2k+1)!!$, we get
\be\label{ujump}
[[u]]=-\lambda^{k+1}\frac{k!}{(2k+1)!}+O(\lambda^{k+2}).
\ee
Multiplying \eq{eig-eq} by $e^{\lambda(y+1)/2}$ and integrating by parts repeatedly, we obtain
\be\label{eig-sol}
u(y)=u(1)e^{-imh-\hf{\lambda}(y+1)}+\hf{(-1)^{k+1}}[[u]]
\sum_{l=0}^\infty(-\hf{\lambda})^lR_{k+1}^{-,(-l)}(y),
\ee
where $R_{k+1}^{-,0}(y)=R_{k+1}^-(y)$, and for $l\ge0$,
\be
R_{k+1}^{-,(-l-1)}(y)=\int_{-1}^yR_{k+1}^{-,(-l)}(z)dz.
\ee
Using the formula
\be
\int_{-1}^y\phi_k(x)dx=\frac{\phi_{k+1}(y)-\phi_{k-1}(y)}{2k+1},\ \ \ k\ge1,
\ee
we have
\be
R_{k+1}^{-,(-l)}=\sum_{i=k-l}^{k+l+1}c^l_i\phi_i,\ \ \ 0\le l\le k,
\ee
where $c^l_i$ are constants. In particular,
\be\label{cll}
c^l_{k-l}=(-1)^{l+1}\frac{(2k-2l+1)!!}{(2k+1)!!}.
\ee
We also have $R_{k+1}^{-,(-l)}(1)=0$ for $0\le l\le k$, and
\be
R_{k+1}^{-,(-k-1)}(1)=2c^k_0=(-1)^{k+1}\frac{2}{(2k+1)!!}.
\ee
Substituting in $y=1$ and using the formula above, we get
\be
u(1)(1-e^{-imh-\lambda})=\hf{(-1)^{k+1}}[[u]]
\left((-\hf{\lambda})^{k+1}(-1)^{k+1}\frac{2}{(2k+1)!!}+O(\lambda^{k+2})\right).
\ee
Substituting in \eq{ujump} we get
\be
\lambda+imh=-\frac{(k+1)!k!}{(2k+2)!(2k+1)!}(mh)^{2k+2}+O((mh)^{2k+3}).
\ee
Substituting the equation above, \eq{ujump}, and $u(1)=e^{imh/2}$ into \eq{eig-sol}, we get
\bes
u(y)&=u(1)e^{-imh+\hf{imh}(y+1)}+\hf{(-1)^{k+1}}[[u]]
\sum_{l=0}^k(-\hf{\lambda})^lR_{k+1}^{-,(-l)}(y)+O((mh)^{2k+2})\\
&=e^{im(x-x_j)}+\hf{(-1)^k}(-imh)^{k+1}\frac{k!}{(2k+1)!}
\sum_{l=0}^k(-\hf{\lambda})^lR_{k+1}^{-,(-l)}(y)+O((mh)^{2k+2})
\ees
By \eq{cll} we see that for $0\le n\le k$,
\be
u_n=p_n^{(m)}+(-1)^{k-n}\frac{(k+1)!k!(2n+1)!}{(2k+2)!(2k+1)!n!}(imh)^{2k+1-n}+O((mh)^{2k+2-n}).
\ee

\end{proof}

Next we prove \thm{Thm-elim} using \lem{Lem-eig}. If the initial error in $u_h(x,0)$ is given by \eq{einit2},
\bes\label{e-dir}
e_j(t)&=\sum_{m=0}^{N-1}e^{imx_j}(e^{-imt}p^{(m)}-e^{A_mt}p^{(m)})\hat{g}_m+O((mh)^{2k+2-n}\hat{g}_m)\\
&=\sum_{m=0}^{N-1}e^{imx_j}\left[(e^{-imt}-e^{\lambda_0^{(m)}t})\tilde{r}_0^{(m)}
-e^{A_mt}\delta^{(m)}+e^{-imt}\delta^{(m)}\right]\hat{g}_m+O((mh)^{2k+2-n}\hat{g}_m).
\ees
Comparing \eqs{xi0}{delta}, we see
\be
-\delta_n^{(m)}\hat{g}_m=\hat{\xi}_{m,n}(0)+O((mh)^{2k+2-n}\hat{g}_m).
\ee
Comparing \eqs{eta-t}{delta}, we see
\be
e^{-imt}\delta^{(m)}\hat{g}_m=\hat{\eta}_{m}(t)+O((mh)^{2k+2-n}\hat{g}_m).
\ee
Lastly, by \eq{lambda0},
\bes\label{er0}
&\sum_{m=0}^{N-1}e^{imx_j}(e^{-imt}-e^{\lambda_0^{(m)}t})\tilde{r}_0^{(m)}\hat{g}_m\\
&=\sum_{m=0}^{N-1}e^{im(x_j-t)}\tilde{r}_0^{(m)}
\left[\frac{(k+1)!k!}{(2k+2)!(2k+1)!}(mh)^{2k+1}mt+O((mh)^{2k+2}mt)\right]\hat{g}_m.
\ees
Comparing \eqs{es}{er0}, we conclude that \eq{e-dir} has the same limits as given in \thm{Thm-elim}. We can see that $\tilde{r}_0^{(m)}$ acts as the special projection $u_I$ for each Fourier mode. In the decomposition,
\be
p^{(m)}=\tilde{c}_0\tilde{r}_0^{(m)}+\sum_{n=1}^k c_n r_n^{(m)},
\ee
by \eq{l0xi}, we have $c_n=O((mh)^{k+1})$ for $1\le n\le k$, and
\be
\tilde{c}_0=1-i(k+1)\frac{(k+1)!k!}{(2k+2)!(2k+1)!}(mh)^{2k+1}+O((mh)^{2k+2}).
\ee
It is interesting that if $\tilde{r}_0^{(m)}$ is replaced by $r_0^{(m)}$ in the decomposition,
\be\label{pr0}
p^{(m)}=c_0r_0^{(m)}+\sum_{n=1}^k c_n r_n^{(m)},
\ee
then $c_0=1+O((mh)^{2k+2})$, due to the following lemma.

\begin{lemma}\label{Lem-r0r0}
The eigenvectors $r_0^{(m)}$ and $\tilde{r}_0^{(m)}$ are related by
\be\label{r0r0}
\tilde{r}_0^{(m)}=r_0^{(m)}\left(1+i(k+1)\frac{(k+1)!k!}{(2k+2)!(2k+1)!}(mh)^{2k+1}+O((mh)^{2k+2})\right).
\ee
\end{lemma}
\begin{proof}
Let
\be
v=\tilde{r}_0^{(m)}\left(1-i(k+1)\frac{(k+1)!k!}{(2k+2)!(2k+1)!}(mh)^{2k+1}\right).
\ee
By \eq{delta},
\be
p_n^{(m)}-v_n=\delta_n^{(m)}
=(-1)^{k+1-n}\frac{(k+1)!k!(2n+1)!}{(2k+2)!(2k+1)!n!}(imh)^{2k+1-n}+O((mh)^{2k+2-n})
\ee
for $1\le n\le k$, and
\be
p_0^{(m)}-v_0
=ik\frac{(k+1)!k!}{(2k+2)!(2k+1)!}(mh)^{2k+1}+O((mh)^{2k+2}).
\ee
Since $p_n^{(m)}$ is real for even $n$ and imaginary for odd n,
\be
\sum_{n=0}^k\frac{(-1)^n(p^{(m)}_n)^2}{2n+1}=\sum_{n=0}^k\frac{|p^{(m)}_n|^2}{2n+1}
=1-\sum_{n=k+1}^\infty\frac{|p^{(m)}_n|^2}{2n+1}=1+O((mh)^{2k+2}).
\ee
Using \eq{pn}, we obtain
\bes
\sum_{n=0}^k\frac{(-1)^nv_n^2}{2n+1}&=v_0^2+\sum_{n=1}^k\frac{(-1)^n(p^{(m)}_n-\delta^{(m)}_n)^2}{2n+1}\\
&=\sum_{n=0}^k\frac{(-1)^n(p^{(m)}_n)^2}{2n+1}
-2\left(p_0^{(m)}(p_0^{(m)}-v_0)+\sum_{n=1}^k(-1)^np_n^{(m)}\delta_n^{(m)}\right)+O((mh)^{2k+2})\\
&=1+O((mh)^{2k+2}).
\ees
By \eq{lr}, $r_0^{(m)}$ satisfies
\be
\sum_{n=0}^k\frac{(-1)^n(r_0^{(m)})_n^2}{2n+1}=1.
\ee
Therefore $\|v-r_0^{(m)}\|=O((mh)^{2k+2})$.
\end{proof}

\subsection{Initialization by special projections}

For $k=1$, as indicated by in \thm{Thm-ebd}, the error in cell average is of order 3 for any $t\in[0,T]$, if $u_h$ is initialized as the $L_2$ projection of $u(x,0)$ onto $V_h$. For $k\ge2$, the error in cell average if of order $k+2$ for small $t$. However, for any interval $[T_0,T]$ with $T_0>0$, the error in cell average for $t\in[T_0,T]$ is of order $2k+1$ for sufficiently small $h$.

\begin{corollary}\label{cor-trans}
Suppose $k\ge2$. For an interval $[T_0,T]$ with $T_0>0$, if
\be\label{hsmall}
h\ln\inv{h}\le \frac{\alpha T_0}{2(k-1)},
\ee
then for any $t\in[T_0,T]$,
\be
\|e_{\cdot,0}(t)\|_{L_2}\lesssim h^{2k+1}(\|g^{(2k+1)}\|_{L_2}+\|g^{(2k+2)}\|_{L_2}t)
\ee
if $g\in H_{2k+2}$, and
\be
\|e_{\cdot,0}(t)\|_\infty\lesssim h^{2k+1}(\|g\|_{2k+1}+\|g\|_{2k+2}t)
\ee
if $g\in H_{2k+3}$.
\end{corollary}

\begin{proof}
For any $t\in[T_0,T]$, \eq{hsmall} implies that
\be
e^{-\frac{\alpha t}{2h}}\le e^{-\frac{\alpha T_0}{2h}}\le h^{k-1}.
\ee
Then the conclusion follows from \thm{Thm-ebd} and \ref{Thm-ebd2}.
\end{proof}

Another way to decrease the exponentially decaying term in the error bound, as shown in Ref. \cite{Cao}, is by initializing $u_h$ to be
\be\label{uIl}
u_I^l=P^-_hu-\sum_{i=1}^lw_i,\ \ \ 0\le l\le k.
\ee

\begin{theorem}\label{Thm-ebdI}
Suppose $g\in H_{2k+3}$, so that $\|g\|_{2k+2}<\infty$. Let $u_h$ be the solution to \eq{numeq} with
\be\label{uinitI}
u_h(x,0)=u_I^l(x,0),
\ee
for some $0\le l\le k$. There exists an $\alpha>0$, such that for any $t\in[0,T]$,
\be
\|\xi(\cdot,t)\|_\infty\lesssim h^{2k+1}(\|g\|_{2k+1}+\|g\|_{2k+2}t)+h^{k+l+2}\|g\|_{k+l+2}
e^{-\frac{\alpha t}{2h}}.
\ee
\be
\|\xi_{\cdot,0}(t)\|_\infty\lesssim h^{2k+1}(\|g\|_{2k+1}+\|g\|_{2k+2}t)+h^{k+l+3}\|g\|_{k+l+3}
e^{-\frac{\alpha t}{2h}}.
\ee
Consequently,
\be
\|e_{\cdot,0}(t)\|_\infty\lesssim h^{2k+1}(\|g\|_{2k+1}+\|g\|_{2k+2}t)+h^{k+l+3}\|g\|_{k+l+3}
e^{-\frac{\alpha t}{2h}},
\ee
\be
\|e^-_\cdot(t)\|_\infty\lesssim h^{2k+1}(\|g\|_{2k+1}+\|g\|_{2k+2}t)+h^{k+l+2}\|g\|_{k+l+2}
e^{-\frac{\alpha t}{2h}},
\ee
and for $1\le n\le k$,
\be
\|e_{\cdot,n}(t)\|_\infty\lesssim h^{2k+1-n}\|g\|_{2k+1-n}+h^{k+l+2}\|g\|_{k+l+2}
e^{-\frac{\alpha t}{2h}}.
\ee
Moreover, for any $t>0$,
\be
\lim_{h\to0}\frac{e_{j,0}(t)}{h^{2k+1}}=\frac{(-1)^k(k+1)!k!}{(2k+2)!(2k+1)!}
[(k-l-1)g^{(2k+1)}(x_j-t)-tg^{(2k+2)}(x_j-t)],
\ee
with uniform convergence, and \eq{en-lim} still holds for $1\le n\le k$.
\end{theorem}

\begin{proof}
\eq{uinitI} implies that $\hat{\xi}_{m,n}(0)$ is the same as given in \eq{xi0} for $0\le n\le k-l-1$, while
\be
\hat{\xi}_{m,n}(0)=O((mh)^{k+l+2}\hat{g}_m),\ \ \ k-l\le n\le k.
\ee
Therefore
\be
l_0^{(m)}\hat{\xi}_m(0)=(k-l)\frac{(k+1)!k!}{(2k+2)!(2k+1)!}
(mh)^{2k+1}i\hat{g}_m+O((mh)^{2k+2}\hat{g}_m),
\ee
and for $1\le n\le k$,
\be
l_n^{(m)}\hat{\xi}_m(0)=O(mh)^{k+l+2}\hat{g}_m).
\ee
The rest of the proof is the same as those for \thm{Thm-ebd} and \ref{Thm-elim}.
\end{proof}

Same result in $L_2$ norm of the errors can be obtained for $g\in H_{2k+2}$. By \thm{Thm-ebdI}, the optimal order of $2k+1$ for error in cell average and $2k+1-n$ in $e_{j,n}$ can be achieved with the initial discretization,
\be\label{init-k2}
u_h(x,0)=P^-_hu(x,0)-\sum_{i=1}^{k-2}w_i(x,0).
\ee
In particular, for $k=2$, we can set $u_h(x,0)=P^-_hu(x,0)$ to achieve the optimal order of 5 in cell average for all $t$. Since the Gauss-Radau projection does not involve time derivative, it is easier to implement than $u^l_I$ with $l>0$. For $k=1$, to achieve optimal superconvergence of order 3 for all $t$, we only need to initialize $u_h$ as the $L_2$ projection of $u(x,0)$, which is even simpler since it is independent of the direction of the flow.

\section{Vector equations}

In this section we study the solution to vector linear advection equations,
\bes
u_t+Au_x=0,\ \ \ (x,t)\in[0,2\pi]\times[0,T].\\
u(x,0)=g(x),\ \ \ u(0,t)=u(2\pi,t).\non
\ees
In the equation, $u$ is a vector, and $A$ is a diagonalizable matrix. If the equation is solved by DG with upwinding flux, it is equivalent to diagonalizing the equation and solving the scalar equation for each eigenmode. To avoid diagonalization, we use the Lax-Friedrichs flux on the cell boundary for the DG method,
\be
F^*(u^-,u^+)=\hf{f(u^-)+f(u^+)}-M\hf{u^+-u^-},
\ee
where $M>0$. For error analysis, it is equivalent to writing $u_h$ as the sum of eigenmodes, each of which satisfying the following equation,
\be\label{numeqM}
((u_h)_t,v)=a(u_h,v_x)+\sum_{j=1}^N([v]F^*(u_h^-,u_h^+)|_{j+\half},\ \ \ \forall v\in V_h,
\ee
where $a$ is the speed of the eigenmode. For $M\neq|a|$, $F^*$ is not the upwinding flux, and the error analysis in the last section has to be modified. In the following theorem, we show that for nonzero $a$, the superconvergence still holds with the same order. Without loss of generality, we set $a=1$ in \eq{numeqM}, and we study the difference between $u_h$ and the exact solution of \eq{master}.

\begin{theorem}\label{Thm-eM}
Suppose $g\in H_{2k+3}$. Let $u_h$ be the solution to \eq{numeqM} with initial error
\be
\hat{e}_{m,n}(0)=O((mh)^{2k+2-n}\hat{g}_m),\ \ \ 0\le n\le k.
\ee
There exists an $\alpha_M>0$, such that for any $t\in[0,T]$,
\be\label{e0M}
\|e_{\cdot,0}(t)\|_\infty\lesssim h^{2k+1}(\|g\|_{2k+1}+\|g\|_{2k+2}t)+h^{k+2}\|g\|_{k+2}
e^{-\frac{\alpha_M t}{2h}},
\ee
and for $1\le n\le k$,
\be\label{enM}
\|e_{\cdot,n}(t)\|_\infty\lesssim h^{2k+1-n}\|g\|_{2k+1-n}+h^{k+1}\|g\|_{k+1}
e^{-\frac{\alpha_M t}{2h}}.
\ee
For any $t>0$,
\be\label{e0-limM}
\lim_{h\to0}\frac{e_{j,0}(t)}{h^{2k+1}}=\chi_M\frac{(-1)^k(k+1)!k!}{(2k+2)!(2k+1)!}
[kg^{(2k+1)}(x_j-t)-tg^{(2k+2)}(x_j-t)],
\ee
and for $1\le n\le k$,
\be\label{en-limM}
\lim_{h\to0}\frac{e_{j,n}(t)}{h^{2k+1-n}}=\chi_M
\frac{(-1)^{k+1-n}(k+1)!k!(2n+1)!}{(2k+2)!(2k+1)!n!}g^{(2k+1-n)}(x_j-t),
\ee
where $\chi_M=M$ for even $k$, and $\chi_M=1/M$ for odd $k$. The convergence is uniform.
\end{theorem}

\begin{proof}
The proof is similar to that of \lem{Lem-eig}. For each Fourier mode and $v\in V_h(\tau_j)$, \eq{numeqM} becomes
\be
((u_h)_t+(u_h)_x,v)=(u_L-u_Re^{-imh})\left(\hf{M-1}v_Re^{imh}-\hf{M+1}v_L\right),
\ee
where $L$ and $R$ represents the left and right end points of $\tau_j$. Let $(u_h)_t=(\lambda/h)u_h$, we can decompose $u_h$ into $k+1$ eigenmodes. Scaled from $x\in\tau_j$ to $y=2(x-x_j)/h\in[-1,1]$, the eigenfunction satisfies
\be\label{eig-eqM}
\lambda u+2u_y=(-1)^{k+1}[[u]](R_{k+1}^M)'(y),\ \ \ u\in P_k.
\ee
where $[[u]]=u(-1)-u(1)e^{-imh}$, and
\be
R^M_{k+1}(y)=\bdf a_M\phi_{k+1}(y)-b_M\phi_k(y),\ \ \ k\ {\rm even}\\
b_M\phi_{k+1}(y)-a_M\phi_k(y),\ \ \ k\ {\rm odd}\edf,
\ee
in which
\be
a_M=e^{\hf{imh}}\left(\cos\hf{mh}-iM\sin\hf{mh}\right),\ \ \
b_M=e^{\hf{imh}}\left(M\cos\hf{mh}-i\sin\hf{mh}\right).
\ee
Notice that $a_M+b_M=M+1$, and $a_M-b_M=(1-M)e^{imh}$. We seek the physical eigenvalue and eigenfunction. Since $u\in P_k$,
\be
u(y)=[[u]]\hf{(-1)^k}\sum_{l=1}^{k+1}\frac{R^{M,(l)}_{k+1}(y)}{(-\hf{\lambda})^l}.
\ee
Substituting in $y=1$ and using the formula $\phi^{(k+1)}_{k+1}=(2k+1)!!$, we get for even $k$,
\be
[[u]]=-\lambda^{k+1}\frac{k!}{(2k+1)!}\frac{u(1)}{a_M}+O(\lambda^{k+2})
=-\lambda^{k+1}\frac{k!}{(2k+1)!}+O(\lambda^{k+2}),
\ee
and for odd $k$,
\be
[[u]]=-\lambda^{k+1}\frac{k!}{(2k+1)!}\frac{u(1)}{b_M}+O(\lambda^{k+2})
=-\lambda^{k+1}\frac{k!}{(2k+1)!}\inv{M}+O(\lambda^{k+2}).
\ee
Multiplying \eq{eig-eqM} by $e^{\lambda(y+1)/2}$ and integrating by parts repeatedly, we obtain
\be
u(y)=\left(u(1)-\hf{M-1}[[u]]e^{imh}\right)e^{-imh-\hf{\lambda}(y+1)}+\hf{(-1)^{k+1}}[[u]]
\sum_{l=0}^\infty(-\hf{\lambda})^lR_{k+1}^{M,(-l)}(y),
\ee
where $R_{k+1}^{M,0}(y)=R_{k+1}^M(y)$, and for $l\ge0$,
\be
R_{k+1}^{M,(-l-1)}(y)=\int_{-1}^yR_{k+1}^{M,(-l)}(z)dz.
\ee
Similar to the proof of \lem{Lem-eig}, we have
\be
R_{k+1}^{M,(-l)}=\sum_{i=k-l}^{k+l+1}c^l_i\phi_i,\ \ \ 0\le l\le k,
\ee
where
\be
c^l_{k-l}=\bdf (-1)^{l+1}\frac{(2k-2l+1)!!}{(2k+1)!!}b_M,\ \ \ k\ {\rm even}\\
(-1)^{l+1}\frac{(2k-2l+1)!!}{(2k+1)!!}a_M,\ \ \ k\ {\rm odd}\edf.
\ee
At $y=1$,
\be
R_{k+1}^{M,(-l)}(1)=\bdf (-1)^k(1-M)e^{imh},&l=0\\0,&1\le l\le k\\2c^k_0,&l=k+1\edf,
\ee
so
\be
\left(u(1)-\hf{M-1}[[u]]e^{imh}\right)(1-e^{-imh-\lambda})=\hf{(-1)^{k+1}}[[u]]
\left((-\hf{\lambda})^{k+1}2c^k_0+O(\lambda^{k+2})\right).
\ee
Substituting in $[[u]]$, we get
\be
\lambda+imh=-\chi_M\frac{(k+1)!k!}{(2k+2)!(2k+1)!}(mh)^{2k+2}+O((mh)^{2k+3}),
\ee
where $\chi_M=M$ for even $k$, and $\chi_M=1/M$ for odd $k$. Normalizing $u$ by
\be\label{normM}
u(1)-\hf{M-1}[[u]]e^{imh}=e^{\hf{imh}},
\ee
we get
\be
u(y)=e^{im(x-x_j)}+\hf{(-1)^{k+1}}[[u]]
\sum_{l=0}^k(-\hf{\lambda})^lR_{k+1}^{M,(-l)}(y)+O((mh)^{2k+2}),
\ee
and so
\be
u_n=p_n^{(m)}+(-1)^{k-n}\chi_M\frac{(k+1)!k!(2n+1)!}{(2k+2)!(2k+1)!n!}(imh)^{2k+1-n}+O((mh)^{2k+2-n}),
\ee
where $p_n^{(m)}$ is the projection of $e^{im(x-x_j)}$ onto $\phi_n(\tau_j)$. The rest of the proof is similar to the argument following the proof of \lem{Lem-eig}. We only need to point out that the nonphysical eigenvalues have negative real parts, because \eq{numeqM} can be written as
\be\label{intpartM}
(u_t+u_x,v)=-\sum_{j=1}^N\left.\left([u](\hf{v^++v^-}+\hf{M}[v])\right)\right|_{j+\half},
\ \ \ \forall v\in V_h.
\ee
It gives the energy estimate,
\be
\diff{}{t}(u,\bar{u})=-M\sum_{j=1}^N|[u]_{j+\half}|^2\le0.
\ee
Since $M>0$, $\Re\lambda^{(m)}_n\le0$. The equality holds only if $[u]=0$, but then $u_t+u_x=0$ by \eq{intpartM}. Since the eigenfunction $u$ is a polynomial, it must be a constant, and $\lambda=0$, which is the physical eigenvalue for $m=0$. Therefore nonphysical eigenvalues have negative real parts. We can let $\alpha_M=\min_{1\le n\le k}\Re(-\lambda_n^{(0)})$.
\end{proof}

If $g\in H_{2k+2}$, the error bounds are given in $L_2$ norm rather than $L_\infty$ norm, and the convergence to the asymptotic error is in $L_2$ norm rather than uniform. \eq{normM} indicates that the downwind error is only of order $k+1$. Similar to \lem{Lem-r0r0}, if the physical eigenvector is normalized by \eq{lr}, the corresponding coefficient in the eigenvector decomposition of $p^{(m)}$, as in \eq{pr0}, is $c_0=1+O((mh)^{2k+2})$.

For $a\neq0$ in \eq{numeqM}, not necessarily 1, the asymptotic error is
\be\label{e0-limMa}
\lim_{h\to0}\frac{e_{j,0}(t)}{h^{2k+1}}={\rm sign}(a)\chi_M\frac{(-1)^k(k+1)!k!}{(2k+2)!(2k+1)!}
[kg^{(2k+1)}(x_j-at)-atg^{(2k+2)}(x_j-at)].
\ee
where $\chi_M=M/|a|$ for even $k$, and $\chi_M=|a|/M$ for odd $k$.

\section{Numerical simulations}

In this section we will perform several numerical experiments to demonstrate the superconvergence properties stated in the previous section. In all examples, the time integration is done by a 5th order Runge-Kutta scheme with the CFL number $0.1$. In the first example, we validate the asymptotic error with smooth initial data, as given in \thm{Thm-elim}.

\begin{exmp}\label{ex1}
\be\begin{aligned}
&u_t(x,t)+u_x(x,t)=0,\ \ \ (x,t)\in[0,2\pi]\times(0,1],\\
&u(x,0)=\sin^{2k+2}(x),\ \ \ u(0,t)=u(2\pi,t).
\end{aligned}\ee
\end{exmp}

\begin{figure}[ht]
\begin{centering}
\begin{subfigure}[b]{0.49\textwidth}
  \includegraphics[width=\textwidth]{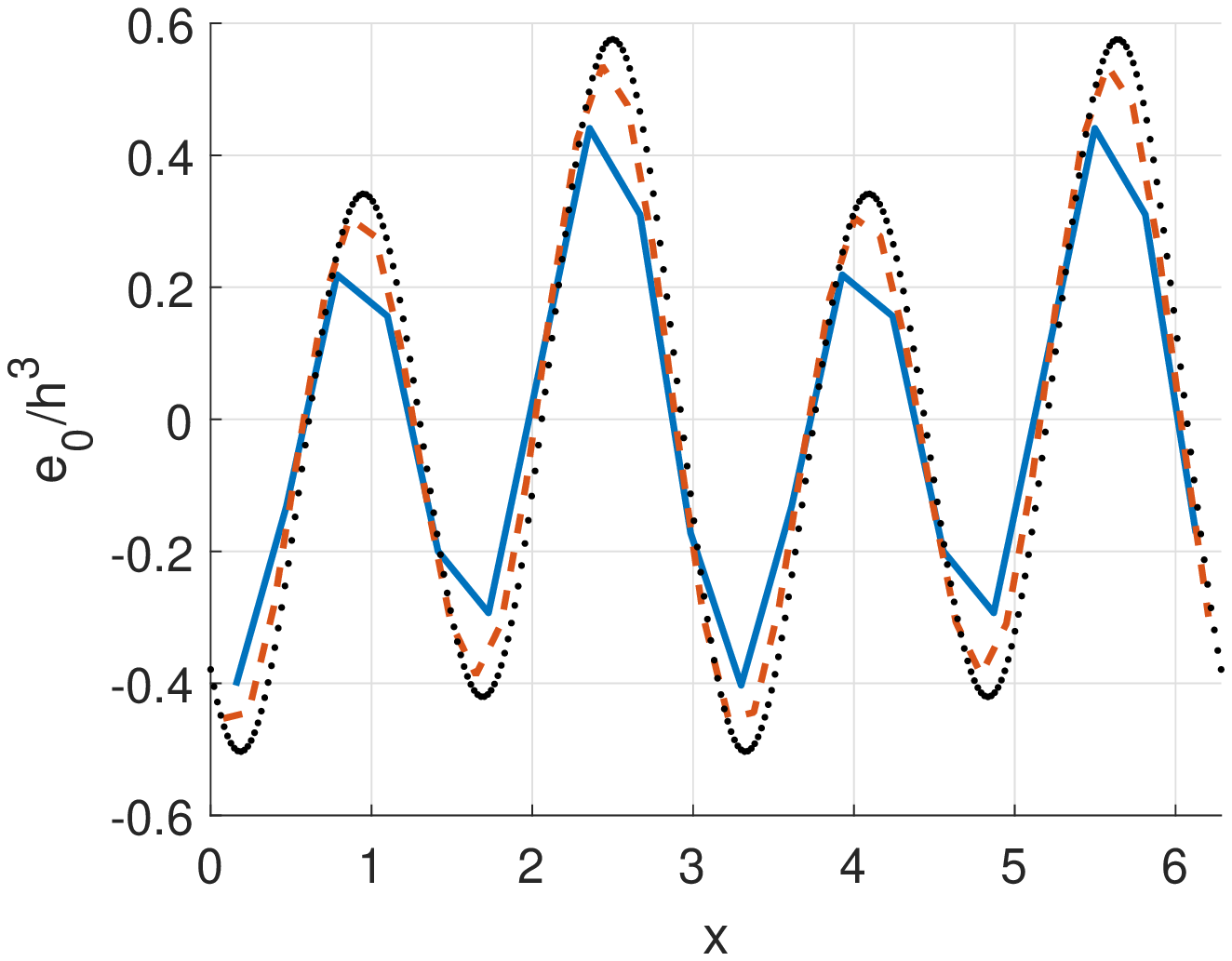}
  \caption{$k=1$}
\end{subfigure}
\begin{subfigure}[b]{0.49\textwidth}
  \includegraphics[width=\textwidth]{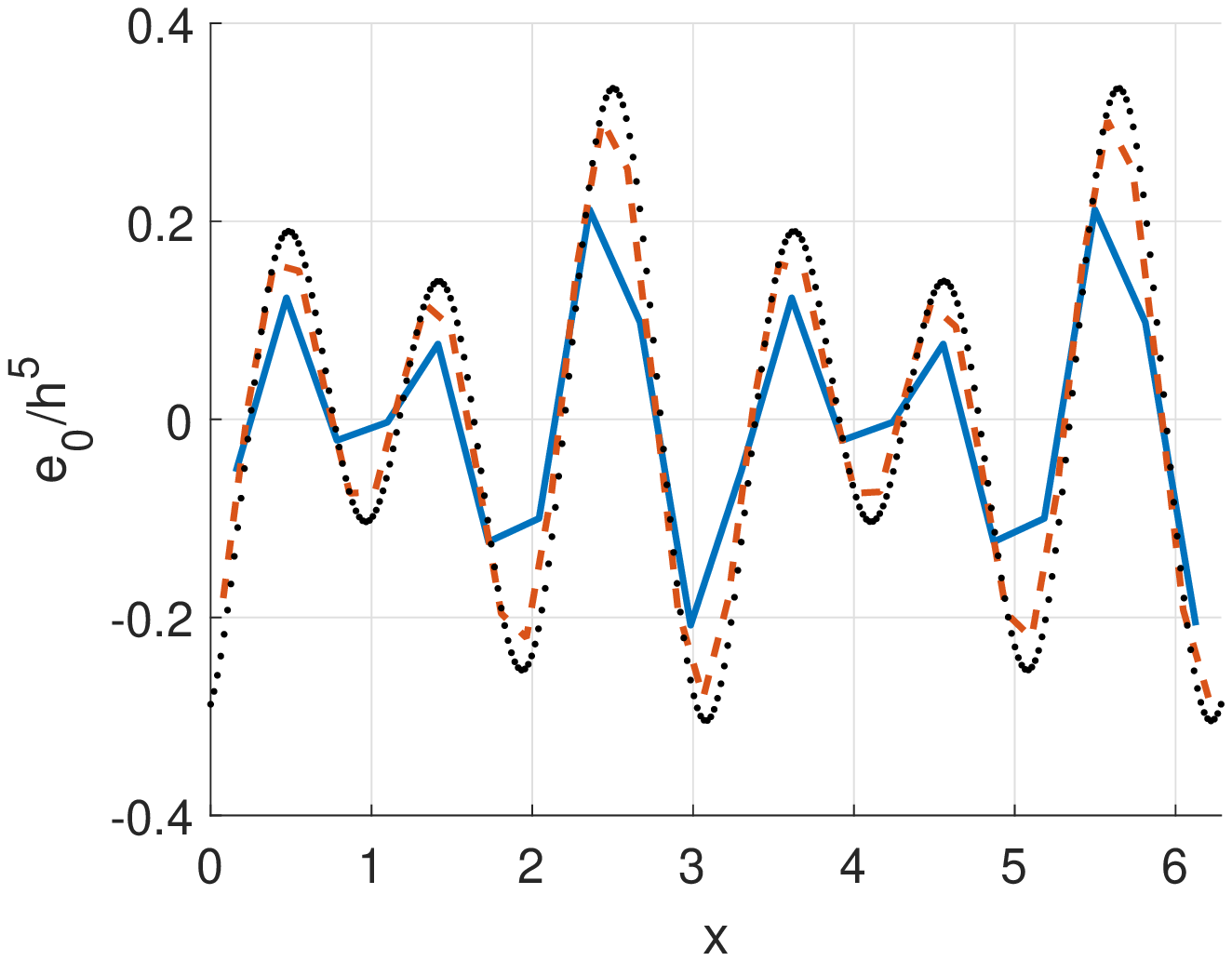}
  \caption{$k=2$}
\end{subfigure}
\end{centering}
\caption{Example 1. The figure shows $e_0/h^{2k+1}$ at $t=1$ for $k=1$ and $k=2$. Solid lines: $N=20$; dashed lines: $N=40$; dotted line: asymptotic error.}\label{fig:ex1}
\end{figure}

\begin{table}[ht]
\caption{Cell average errors in Example 1 for $k=1,2,3$.}\label{tab:ex1}
\centering
\begin{tabular}{|c |c ||c|c|| c|c|| c|c|}
\hline
 $k$  & N & $\|e_0\|_1$ & order & $\|e_0\|_2$ & order& $\|e_0\|_\infty$ & order\\
\hline
  &  40&1.10E-03 & -    & 1.20E-03 & -    & 2.10E-03 & -   \\
1&  80& 1.41E-04 & 2.93 & 1.59E-04 & 2.93 & 2.73E-04 & 2.93\\
  &160& 1.79E-05 & 2.98 & 2.02E-05 & 2.98 & 3.47E-05 & 2.98\\
  &320& 2.25E-06 & 2.99 & 2.54E-06 & 2.99 & 4.35E-06 & 2.99\\
\hline
  & 40& 1.28E-04 & -    & 1.52E-04 & -    & 2.87E-04 & -   \\
2&  80& 4.24E-07 & 4.91 & 5.02E-07 & 4.92 & 9.68E-07 & 4.89\\
  &160& 1.35E-08 & 4.98 & 1.59E-08 & 4.98 & 3.07E-08 & 4.98\\
  &320& 4.19E-10 & 5.01 & 4.94E-10 & 5.01 & 9.55E-10 & 5.01\\
\hline
  & 40& 2.60E-07 & -    & 3.24E-07 & -    & 6.90E-07 & -   \\
3&  80& 1.32E-09 & 7.62 & 1.68E-09 & 7.59 & 3.68E-09 & 7.55\\
  &160& 1.05E-11 & 6.98 & 1.29E-11 & 7.02 & 2.69E-11 & 7.09\\
  &320& 7.65E-14 & 7.10 & 9.43E-14 & 7.10 & 1.95E-13 & 7.11\\
\hline
\end{tabular}
\end{table}

\fig{fig:ex1} plots the numerical error in cell average at $t=1$ for \exm{ex1} with $k=1$ and $k=2$, along with the asymptotic errors given in \thm{Thm-elim}. The initial discretization is the $L_2$ projection of $u(x,0)$. The figure shows that the cell average error is of order $2k+1$ as $h\to0$, and converges to the asymptotic error given in \eq{e0-lim}. \tab{tab:ex1} show the cell average errors in \exm{ex1} for $k=1,2,3$. Since the initial data is in $H^{2k+3}$, the cell average error at $t=1$ in $L_1$, $L_2$, and $L_\infty$ norms are all of order $2k+1$ as $h\to0$.

The next example shows that for non-smooth initial data, the error may have lower order.

\begin{exmp}\label{ex2}
\be\begin{aligned}
&u_t(x,t)+u_x(x,t)=0,\ \ \ (x,t)\in[0,2\pi]\times(0,1],\\
&u(x,0)=|\sin^{2k+1}(x)|,\ \ \ u(0,t)=u(2\pi,t).
\end{aligned}\ee
\end{exmp}

\begin{figure}[ht]
\begin{centering}
\begin{subfigure}{0.49\textwidth}
  \includegraphics[width=\textwidth]{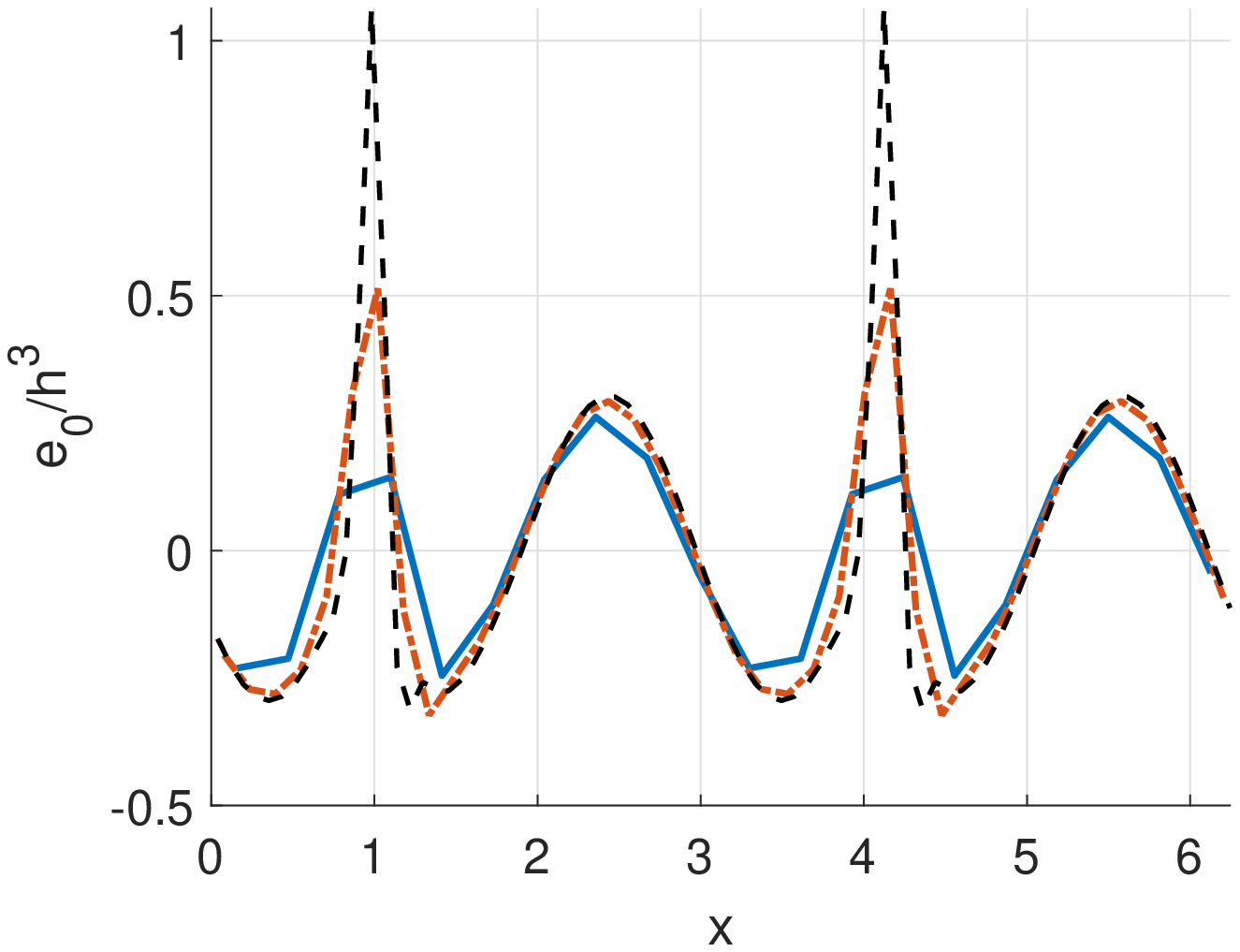}
  \caption{$k=1$}
\end{subfigure}
\begin{subfigure}{0.49\textwidth}
  \includegraphics[width=\textwidth]{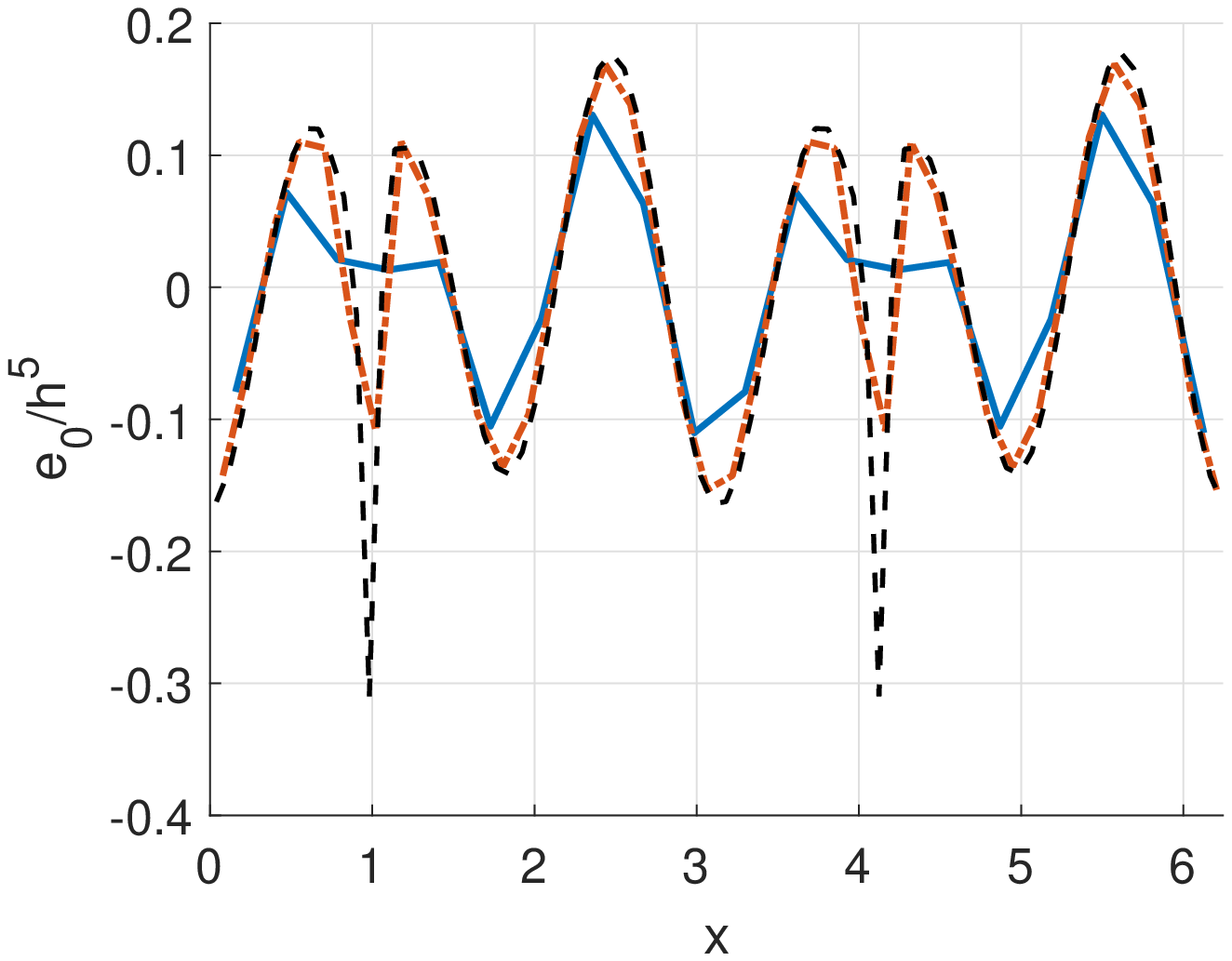}
  \caption{$k=2$}
\end{subfigure}
\end{centering}
\caption{Example 2. The figure shows $e_0/h^{2k+1}$ at $t=1$ for $k=1$ and $k=2$. Solid lines: $N=20$; dash-dotted lines: $N=40$; dashed line: $N=80$.}\label{fig:ex2}
\end{figure}

In this example, $u(x,0)\in H^{2k+1}$ but not in $H^{2k+2}$. As shown in \fig{fig:ex2}, the cell average errors are of order $2k$ at points where $u^{2k+1}(x,t)$ is discontinuous, and of order $2k+1$ elsewhere. As a result, the $L_\infty$ norm of the cell average error is of order $2k$. Due to the localness of the spikes in \fig{fig:ex2}, $\|e_0\|_2$ is of order $2k+1/2$, and $\|e_0\|_1$ is still of order $2k+1$, as demonstrated in \tab{tab:ex2}. The error $\|e_0\|_2$ appears to have order higher than $2k+1/2$ in \tab{tab:ex2} because the error near the spikes dominates over the error elsewhere only for small $h$.

\begin{table}
\caption{Cell average errors in Example 2 for $k=1,2,3$.}\label{tab:ex2}
\centering
\begin{tabular}{|c |c ||c|c|| c|c|| c|c|}
\hline
 $k$  & N & $\|e_0\|_1$ & order & $\|e_0\|_2$ & order& $\|e_0\|_\infty$ & order\\
\hline
  & 40& 8.14E-04 & -    & 9.22E-04 & -    & 1.99E-03 & -   \\
1&  80& 1.12E-04 & 2.86 & 1.39E-04 & 2.73 & 5.15E-04 & 1.95\\
  &160& 1.46E-05 & 2.94 & 2.10E-05 & 2.73 & 1.14E-04 & 2.17\\
  &320& 1.86E-06 & 2.98 & 3.21E-06 & 2.71 & 2.41E-05 & 2.25\\
\hline
  & 40& 8.66E-06 & -    & 9.75E-06 & -    & 1.62E-04 & -   \\
2&  80& 2.93E-07 & 4.88 & 3.42E-07 & 4.83 & 9.26E-07 & 4.13\\
  &160& 9.58E-09 & 4.94 & 1.18E-08 & 4.86 & 5.49E-08 & 4.08\\
  &320& 3.03E-10 & 4.98 & 4.16E-10 & 4.83 & 3.02E-09 & 4.18\\
\hline
  & 40& 1.71E-07 & -    & 2.09E-07 & -    & 4.13E-07 & -   \\
3&  80& 9.18E-10 & 7.54 & 1.07E-09 & 7.61 & 2.12E-09 & 7.61\\
  &160& 6.93E-12 & 7.05 & 8.12E-12 & 7.04 & 2.18E-11 & 6.61\\
  &320& 5.04E-14 & 7.10 & 6.14E-14 & 7.05 & 2.98E-13 & 6.19\\
\hline
\end{tabular}
\end{table}

We investigate the effect of initial discretization in the following example.

\begin{exmp}\label{ex3}
\be\begin{aligned}
&u_t(x,t)+u_x(x,t)=0,\ \ \ (x,t)\in[0,2\pi]\times(0,1],\\
&u(x,0)=\left(\frac{x(2\pi-x)}{2\pi}\right)^{2k+2},\ \ \ u(0,t)=u(2\pi,t).
\end{aligned}\ee
\end{exmp}

\begin{figure}[ht]
\begin{centering}
\begin{subfigure}{0.32\textwidth}
  \includegraphics[width=\textwidth]{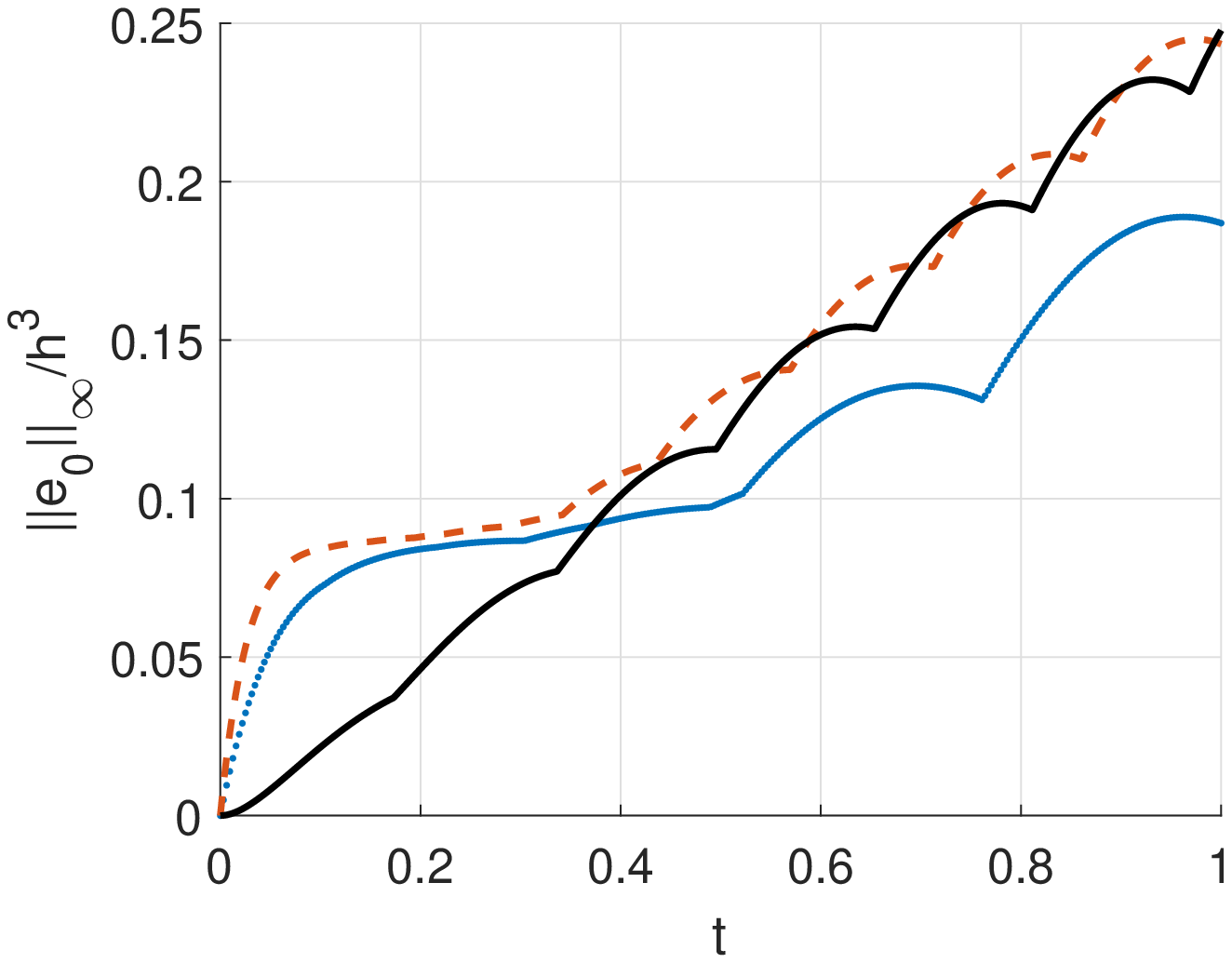}
  \caption{$k=1$}
\end{subfigure}
\begin{subfigure}{0.32\textwidth}
  \includegraphics[width=\textwidth]{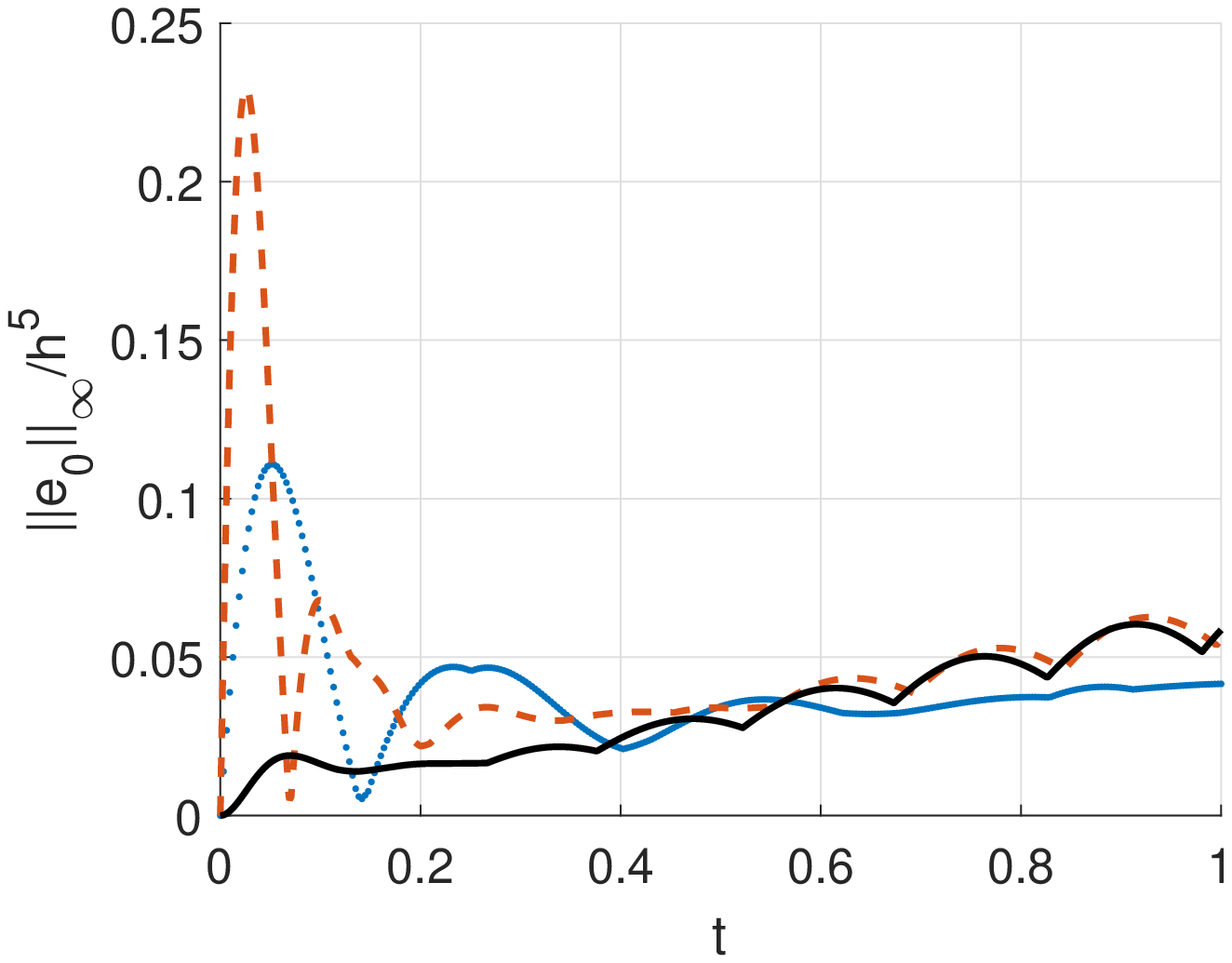}
  \caption{$k=2$}
\end{subfigure}
\begin{subfigure}{0.32\textwidth}
  \includegraphics[width=\textwidth]{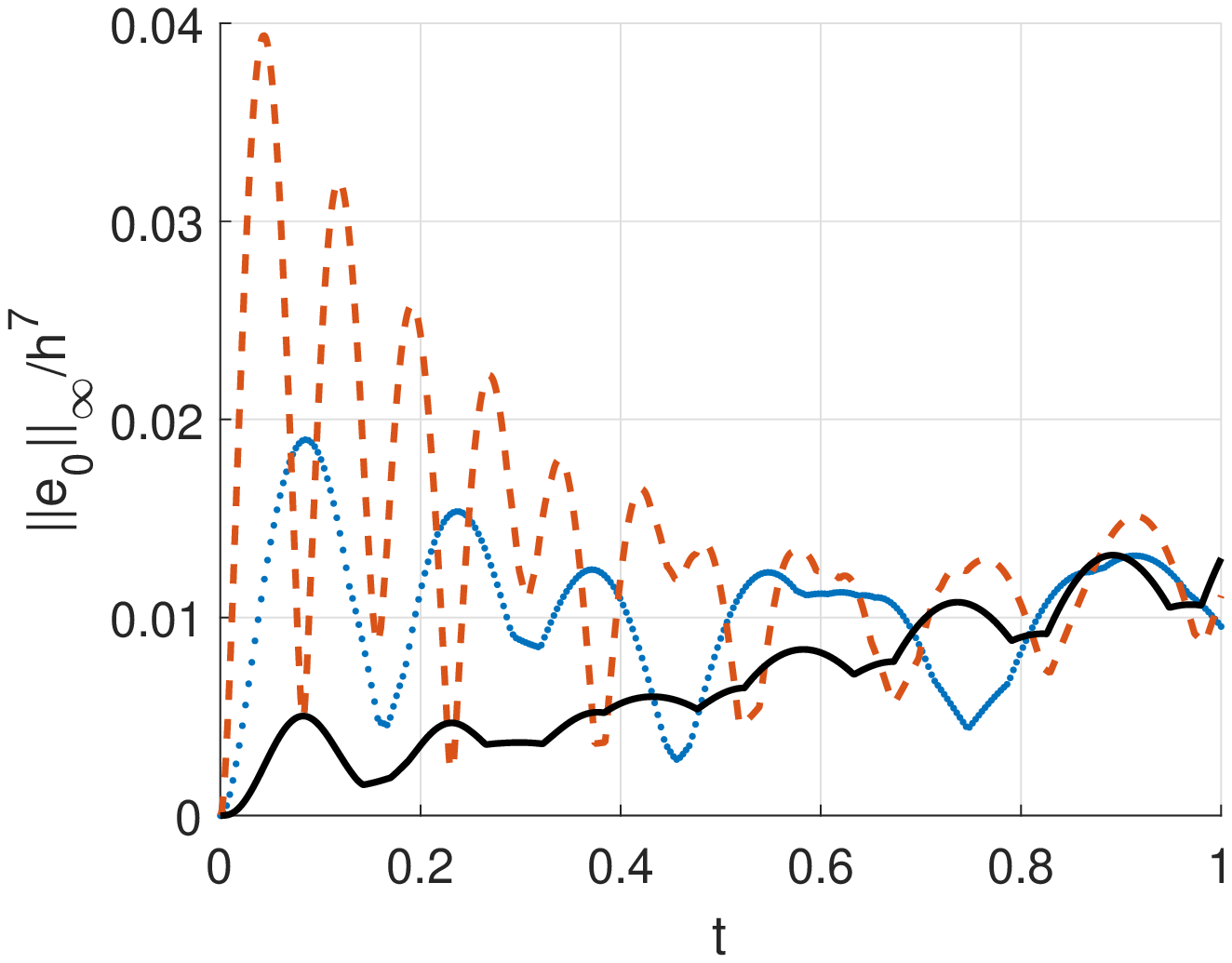}
  \caption{$k=3$}
\end{subfigure}
\end{centering}
\caption{Example 3. The figure shows $\|e_0\|_\infty/h^{2k+1}$ for $k=1,2,3$. For $k=1,2$, solid lines: $N=40$ with Gauss-Radau initialization; dashed(dotted) lines: $N=40(20)$ with $L_2$ initialization. For $k=3$, solid line: $N=40$ with $u_h(x,0)=u_I^1(x,0)$; dashed(dotted) line: $N=40(20)$ with Gauss-Radau initialization.}\label{fig:ex3}
\end{figure}

The initial data is chosen to be in $H^{2k+3}$ and contains infinitely many Fourier modes. The cell average error at any $t>0$ has order $2k+1$ as $h\to0$. However, for a fixed $h$, the error at $t=O(h)$ may have a lower order, depending on the initial discretization. \fig{fig:ex3} shows $\|e_0\|_\infty$ as functions of time for $k=1,2,3$. For $k=1$, the error is of order 3 at any $t$, disregard of the initial discretization. For $k=2$, if $L_2$ initialization is used, the error is of order 4 at small $t$; while for Gauss-Radau initialization, the error if of order 5 at any $t$. For $k=3$, the error is of order 6 at small $t$ for Gauss-Radau initialization, and of order 7 at any $t$ for $u_h(x,0)=u_I^1(x,0)$ as defined in \eq{uIl}. In \fig{fig:ex3}, the transient error for $k=3$ lasts for much longer time than that for $k=2$, because $\alpha=3$ for $k=2$, while $\alpha=0.42$ for $k=3$.

Next example is the linearized Euler equations for isothermal gas. It demonstrates the superconvergence of numerical errors for vector equations solved using the Lax-Friedrichs flux.

\begin{figure}[ht]
\begin{subfigure}{0.49\textwidth}
  \includegraphics[width=\textwidth]{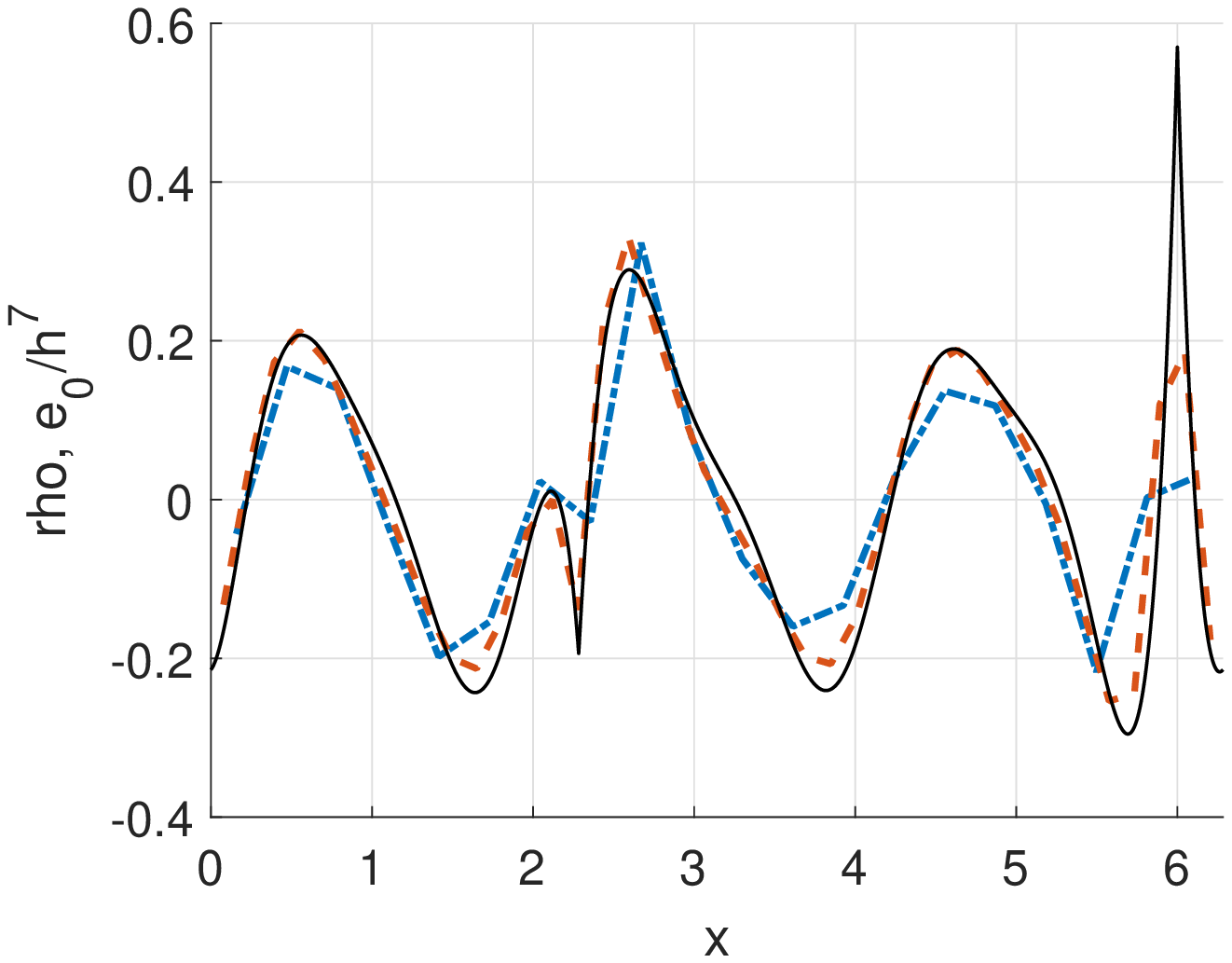}
  \caption{$e_0/h^7$ for $\rho$}
\end{subfigure}
\begin{subfigure}{0.49\textwidth}
  \includegraphics[width=\textwidth]{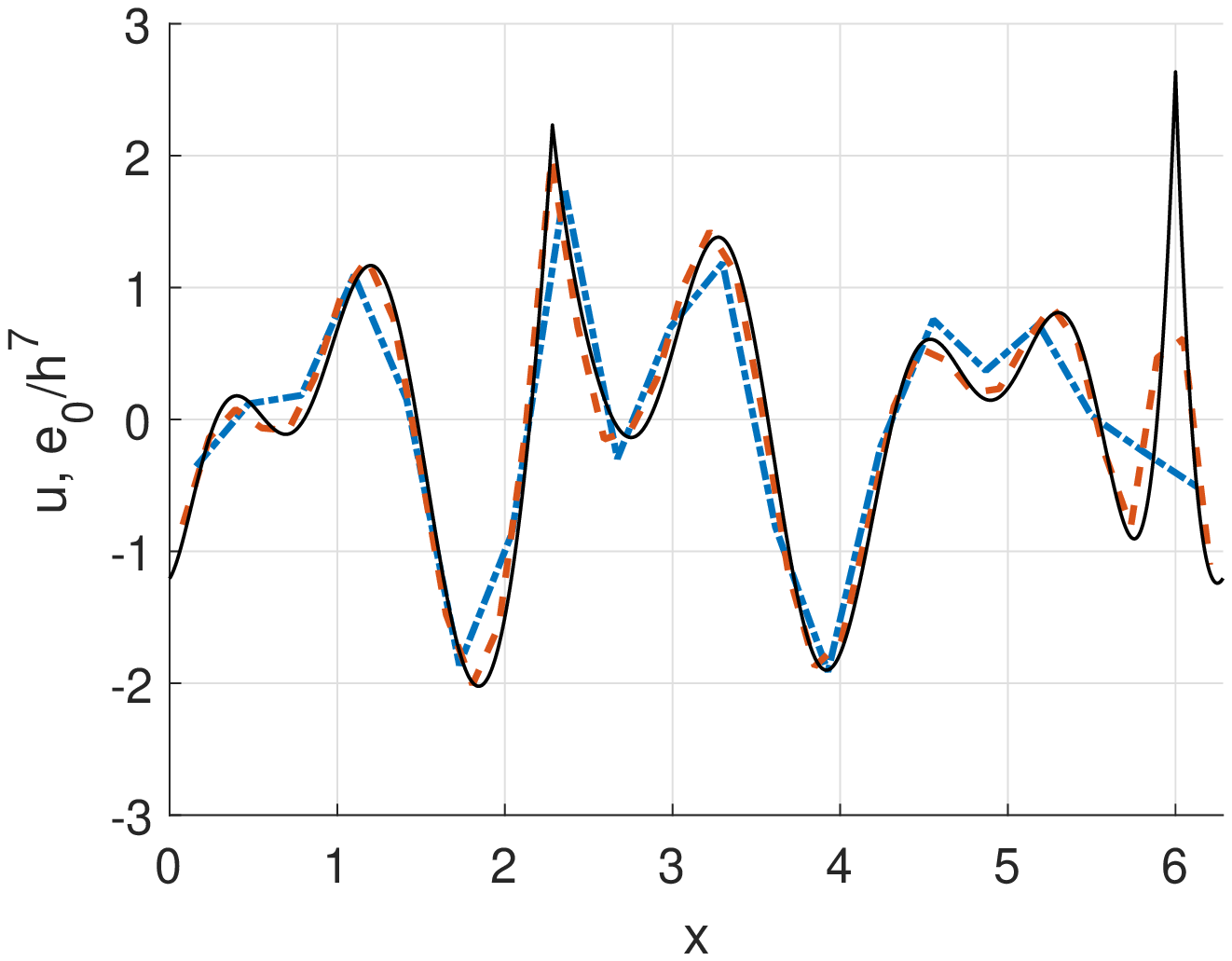}
  \caption{$e_0/h^7$ for $u$}
\end{subfigure}
\caption{Example 4. The cell average error of $\rho$ and $u$ at $t=1$ for $k=3$. Dash-dotted lines: $N=20$; dashed lines: $N=40$; solid line: asymptotic error.}\label{fig:ex4}
\end{figure}

\begin{exmp}\label{ex4}
\be\begin{aligned}
&\rho_t+u_0\rho_x+\rho_0u_x=0,\ \ \rho_0(u_t+u_0u_x)+c^2\rho_x=0,
\ \ \ (x,t)\in[0,2\pi]\times(0,1];\\
&\rho(x,0)=\sin^6(x),\ \ u(x,0)=\left(\frac{x(2\pi-x)}{4}\right)^8;
\ \ \ \rho(0,t)=\rho(2\pi,t),\ \ u(0,t)=u(2\pi,t).
\end{aligned}\ee
\end{exmp}
We set $\rho_0=1$, $u_0=1$, $c=5$, and so there two waves moving at speed 6 and $-4$ respectively. The equation is solved using th Lax-Friedrichs flux with $M=6$. The initialization is done by the $L_2$ projection. Unlike the Gauss-Radau projection, diagonalization and different treatment for right and left going waves is not needed. The problem is solved by DG with $k=1,2,3$. The combined $L_2$ cell average error of $\rho$ and $u$ at $t=1$ is listed in \tab{tab:ex4}. It confirms the superconvergence of order $2k+1$. It's interesting to notice that the error for $N=640$ and $k=1$ is close to that for $N=20$ and $k=3$, while the former takes 100 times longer computational time than the latter. For $k=3$, the cell average errors of $\rho$ and $u$ at $t=1$ are plotted in \fig{fig:ex4}. It shows the convergence of the errors to the asymptotic cell average errors, which are computed by applying \eq{e0-limMa} to both left and right going waves.

\begin{table}[ht]
\caption{Cell average errors in Example 4 for $k=1,2,3$.}\label{tab:ex4}
\centering
\begin{tabular}{|c|c|c|c|c|c|c|c|c|}
\hline
\multicolumn{3}{|c|}{$k=1$}&\multicolumn{3}{c}{$k=2$}&\multicolumn{3}{|c|}{$k=3$}\\
\hline
  N& $\|e_0\|_2$ & order & N& $\|e_0\|_2$ & order & N& $\|e_0\|_2$ & order\\
\hline
 80& 7.72E-02& -    & 40 & 1.10E-03 & -    & 20 & 1.96E-04 & -    \\
160& 9.70E-03& 2.99 & 80 & 3.92E-05 & 4.85 & 40 & 1.57E-06 & 6.97 \\
320& 1.20E-03& 3.00 &160 & 1.29E-06 & 4.93 & 80 & 1.22E-08 & 7.01 \\
640& 1.52E-04& 3.00 &320 & 4.15E-08 & 4.95 &160 & 1.03E-10 & 6.88 \\
\hline
\end{tabular}
\end{table}

\section{Conclusion}
In this paper, we studied the superconvergence of the semi-discrete discontinuous Galerkin method for scalar and vector linear advection equations in one spatial dimension. We used Fourier analysis to prove that the numerical error consists of an asymptotic part and a transient part that decay in time exponentially. For the cell average, the asymptotic error grows linearly in time and has order $2k+1$; while the error projected onto the $n$-th order Legendre polynomial has order $2k+1-n$. The order of the transient part depends on the initial discretization: $k+1$ for $L_2$ projection, $k+2$ for Gauss-Radau projection, etc. The transient error of cell average is one order higher, ie., $k+2$ for $L_2$ projection, $k+3$ for Gauss-Radau projection, etc. We derived the asymptotic error in two ways. In the first approach, we solved the equation for the deviation of the numerical solution from a special interpolation of the exact solution. In the second approach, we decomposed the $L_2$ projection of the initial data into physical and non-physical modes, and computed the asymptotic error by analyzing the physical eigenvalue and eigenvector. Both approaches gave the same asymptotic error, which depends on the initial discretization. Then we extended the Fourier analysis to vector advection equations. Lax-Friedrichs flux was used in order to avoid diagonalization. We showed that the error bounds and asymptotic errors are of the same order, but modified depending on the parity of $k$. All the theoretical results presented have been validated by numerical examples.

Although the current work is on linear advection equations with periodic boundary condition solved by DG on a uniform mesh, numerical experiment shows that much of the results can be extended to more general settings. Our future work involves the analysis of superconvergence for nonlinear advection equations with physical boundary conditions solved on nonuniform grids.


\begin{thebibliography}{}
\bibitem{Cao}
W. Cao, C. Shu, Y. Yang, Z. Zhang.
Superconvergence of Discontinuous Galerkin Method for Scalar Nonlinear Hyperbolic Equations.
{\it SIAM J. Numer. Anal.}, {\bf 56(2)}:732-765 (2018).

\bibitem{Cha}
N. Chalmers, L. Krivodonova.
Spatial and Modal Superconvergence of the Discontinuous Galerkin Method for Linear Equations.
{\it J. Sci. Comput.}, {\bf 72}:128-146 (2017).

\bibitem{Reed}
W.H. Reed, T.R. Hill.
Triangular Mesh for Neutron Transport Equation.
\textit{Los Alamos Scientific Laboratory Report}, \textbf{LA-UR}:73-479 (1973).

\bibitem{Cockburn1}
B. Cockburn, C.W. Shu.
TVB Runge-Kutta local projection discontinuous Galerkin finite element method for conservation laws II: General framework.
{\it Math. Comp.}, {\bf 52}:411-435 (1989).

\bibitem{Cockburn2}
B. Cockburn, S.Y. Lin, C. Shu.
TVB Runge-Kutta local projection discontinuous Galerkin finite element method for conservation laws III: One dimensional systems. \textit{J. Comput. Phys.}, \textbf{84}:90-113 (1989).

\bibitem{Cockburn3}
B. Cockburn, S. Hou, C. Shu.
The Runge-Kutta local projection discontinuous Galerkin finite element method for conservation laws IV, The multidimensional case.
\textit{Math. Comp.}, \textbf{54}:545-581 (1990).

\bibitem{Cockburn4}
B. Cockburn, C.W. Shu.
The Runge-Kutta discontinuous Galerkin method for conservation laws V: Multidimensional systems.
\textit{J. Comput. Phys.}, \textbf{141}:199-224 (1998).

\bibitem{Adjerid1}
S. Adjerid, M. Baccouch.
The discontinuous Galerkin method for two-dimensional hyperbolic problems. Part I: superconvergence error analysis.
\textit{J. Sci. Comput.}, \textbf{33}:75-113 (2007).

\bibitem{Adjerid2}
S. Adjerid, M. Baccouch.
The discontinuous Galerkin method for two-dimensional hyperbolic problems Part II: a posteriori error estimation.
\textit{J. Sci. Comput.}, \textbf{38}:15-49 (2009).

\bibitem{Adjerid3}
S. Adjerid, M. Baccouch.
Asymptotically exact a posteriori error estimates for a one-dimensional linear hyperbolic problem.
\textit{Appl. Numer. Math.}, \textbf{60}:903-914 (2010).

\bibitem{Cao1}
W. Cao, Z. Zhang, Q. Zou.
Superconvergence of discontinuous Galerkin method for linear hyperbolic equations.
\textit{SIAM J. Numer. Anal.}, \textbf{5}:2555-2573 (2014).

\bibitem{Cheng}
Y. Cheng, C.W. Shu.
Superconvergence of discontinuous Galerkin and local discontinuous Galerkin	schemes for linear hyperbolic and convection-diffusion equations in one space dimension.
\textit{SIAM J.}, \textbf{47}:4044-4072 (2010).
	
\bibitem{Yang}
Y. Yang, C.W. Shu.
Analysis of optimal supercovergence of discontinuous Galerkin method for linear hyperbolic equations.
\textit{SIAM J. Numer. Anal.}, \textbf{50}:3110-3133 (2012).

\bibitem{Baccouch}
M. Baccouch.
Recovery-based error estimator for the discontinuous Galerkin method for nonlinear scalar conservation laws in one space dimension.
\textit{J. Sci. Comput.}, \textbf{66}:459-476 (2016).
	
\bibitem{Meng}
X. Meng, C.W. Shu, Q. Zhang, B. Wu.
Superconvergence of discontinuous Galerkin methods for scalar nonlinear conservation laws in one space dimension.
\textit{SIAM J.}, \textbf{50}:2336–2356	(2012).

\bibitem{Bramble}
J.H. Bramble, A.H. Schatz.
Higher order local accuracy by averaging in the finite element method.
\textit{Math Comput.}, \textbf{31}:94-111 (1977).
	
\bibitem{Cockburn5}
B. Cockburn, M. Luskin, C.W. Shu, E. Süli.
Enhanced accuracy by post-processing for finite element methods for hyperbolic equations.
\textit{Math. Comput.}, \textbf{72}:577-606 (2003).

\bibitem{Mirzaee}
H. Mirzaee, L. Ji, J.K. Ryan, R.M. Kirby.
Smoothness-Increasing Accuracy-Conserving (SIAC) postprocessing for discontinuous Galerkin solutions over structured triangular meshes.
\textit{SIAM J.}, \textbf{49}:1899-1920 (2011).

\bibitem{Ji}
L. Ji, Y. Xu, J.K. Ryan.
Accuracy-enhancement of discontinuous Galerkin solutions for convectiondiffusion equations in multiple-dimensions.
\textit{Math. Comput.}, \textbf{81}:1929-1950 (2012).
	
\bibitem{Ji1}
L. Ji, Y. Xu, J.K. Ryan.
Negative-order norm estimates for nonlinear hyperbolic conservation laws.
\textit{J. Sci. Comput.}, \textbf{54}:531-548 (2013).
	
\bibitem{Adjarid}
S. Adjerid, T.C. Massey.
Superconvergence of discontinuous Galerkin solutions for a nonlinear scalar hyperbolic problem.
\textit{Comput. Methods Appl. Mech. Engrg.}, \textbf{195}:3331-3346 (2006).
		
\bibitem{Guo}
W. Guo, X. Zhong, J. Qiu.
Superconvergence of discontinuous Galerkin and local discontinuous Galerkin methods: Eigen-structure analysis based on Fourier approach.
\textit{Journal of Computational Physics}, \textbf{235}:458-485 (2013).

\bibitem{Frean}
D. Frean, J. Ryan.
Superconvergence and the numerical flux: A study using the upwind-biased flux in discontinuous Galerkin methods.
\textit{Comm. on Appl. Math. and Comput.}, \textbf{2}:461-486 (2020).
	
\bibitem{Xu}
Y. Xu, X. Meng, C. Shu, Q. Zhang.
Superconvergence Analysis of the Runge–Kutta Discontinuous Galerkin Methods for a Linear Hyperbolic Equation.
\textit{Journal of Scientific Computing}, \textbf{84}:23 (2020).

\bibitem{Kri}
L. Krivodonova, R. Qin.
An analysis of the spectrum of the discontinuous Galerkin method.
{\it Applied Numerical Mathematics}, {\bf 64}:1-18 (2013).

\bibitem{Abr}
M. Abramowitz, I. Stegun, editors.
{\it Handbook of Mathematical Functions}, Dover, New York (1965).

\end{thebibliography}
\end{document}